\documentclass[11pt]{article}

\usepackage{a4wide}
\usepackage{theorem}
\usepackage{amsmath}
\usepackage{array}
\usepackage{amssymb}
\usepackage{amsfonts}
\usepackage[english]{babel}
\usepackage{epsf}
\usepackage{epsfig}
\usepackage{graphicx}

\newcommand{\R}{\mathbb{R}}
\newcommand{\N}{\mathbb{N}}

\newcommand{\E}{\mathbb{E}}
\newcommand{\Z}{\mathbb{Z}}

\newcommand{\Sd}{\textup{Sd}} 
 
\newcommand{\St}{\textup{St}} 
\newcommand{\Lk}{\textup{Lk}} 
\newcommand{\CW}{\textup{CW}} 
\newcommand{\dvol}{d\textup{vol}} 
\newcommand{\Bd}{\partial} 
 
\newcommand{\ind}{\textup{ind}} 
\newcommand{\ini}{\textup{in}} 
\newcommand{\fin}{\textup{fin}}

\newcommand{\sDp}{\stackrel{\circ}{\Delta}_p}

\newcommand{\mC}{\mathcal{C}}
\newcommand{\mmC}{\bold{C}}

\theoremstyle{definition}
\newtheorem{thm}{Theorem}
\newtheorem{lem}[thm]{Lemma}
\newtheorem{prop}[thm]{Proposition}
\newtheorem{cor}[thm]{Corollary}
\newtheorem{rem}[thm]{Remark}
\newtheorem{ex}[thm]{Example}
\newtheorem{defn}[thm]{Definition}
\newenvironment{pf}{{\bf{Proof.}}}{\hfill $\Box$\\}

\begin{document}

\title{Asymptotic topology of random subcomplexes in a finite simplicial  complex}

\author{Nerm{\accent95\i}n Salepc{\accent95\i} and Jean-Yves Welschinger}

\maketitle

\begin{abstract}
We consider a finite simplicial complex  $K$ together with its successive barycentric subdivisions $\Sd^d(K), d\geq0,$ and study the expected topology of a random subcomplex in $\Sd^d(K), d\gg0$. We get asymptotic upper and lower bounds for the expected Betti numbers of  those subcomplexes, together with the average Morse inequalities and expected Euler characteristic. 

\vspace{0.5cm}


\end{abstract}

\section{Introduction}
\begin{figure}[h]
   \begin{center}
    \includegraphics[scale=0.4]{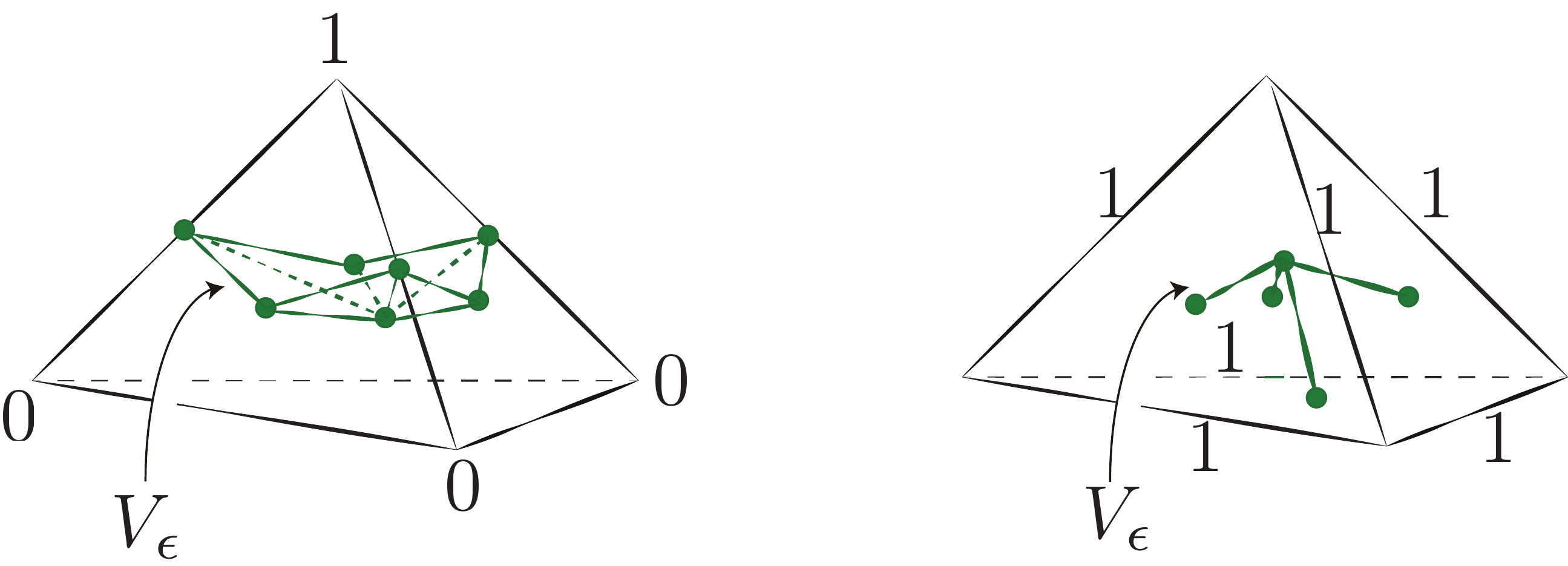}
    \caption{Some examples of simplices of $V_\epsilon$. On the left the case $n=3, k=1$; on the right $n=3, k=2$.}
    \label{Vepsilon}
      \end{center}
 \end{figure}

Let $K$ be a locally finite simplicial complex of dimension $n$ and $\Sd(K)$ be its first barycentric subdivision. Let  $\mC^{k-1}(K)$  be the group of $(k-1)$-dimensional simplicial cochains of $K$ with $\Z/2\Z$-coefficients, $k\in \{1,2,\ldots,n\}$.
For every $\epsilon \in \mC^{k-1}(K)$, we denote by $V_\epsilon$ 
the subcomplex of $\Sd(K)$ 
dual to the cocycle $d \epsilon$, where $d: \mC^{k-1}(K)\to \mC^{k}(K)$ denotes the coboundary operator, see \cite{M}. 
Recall that simplices of  $\Sd(K)$ of dimension $i \in \{0,1,\ldots, n\}$
are of the form $[\hat{\sigma}_0,\ldots, \hat{\sigma}_i]$, where
 $\hat{\sigma}_j$ denotes the barycenter of the simplex $\sigma_j\in K$ for all $j\in \{0,\ldots,i\}$ and $\sigma_j<\sigma_{j+1}$, that is $\sigma_j$ is a proper face of $\sigma_{j+1}$ for all $j\in \{0,\ldots,i-1\}$.
A simplex  thus belongs to  $V_\epsilon$ if and only if it is a face of a simplex $[\hat{\sigma}_0,\ldots, \hat{\sigma}_i]$ such that $\dim\sigma_0=k$ and $<d\epsilon, \sigma_0>\neq 0$. The latter condition means that  $\epsilon$ must take value 1 on an odd number of facets of $\sigma_0$, see Figure~\ref{Vepsilon}.  In other words, $V_\epsilon$ is the union of the blocks $\overline{D}(\sigma_0)$ dual to the simplices $\sigma_0\in K$ such that $<d\epsilon, \sigma_0>\neq 0$, see Section~\ref{Sect_$Vepsilon$} and \cite{M}.

When $k=1$ and $K$ is the moment polytope of some toric manifold equipped with a convex triangulation,  the pair $(K, V_\epsilon)$  gets homeomorphic to the pair  $(K, V'_\epsilon)$, where $V'_\epsilon$  is the  patchwork (tropical) hypersurface defined by O.~Viro, see Proposition~\ref{Prop_$k=1$} and \cite{V2, V1}.

Of special interest are triangulations of compact (topological) manifolds. However, when $k>1$ and $\epsilon\in \mC^{k-1}(K),$
$V_\epsilon$ does not inherit the structure of a triangulated codimension $k$ submanifold, see Remark~\ref{Rem_VHom}. When $k=1$, we prove the following theorem (see Corollary~\ref{Cor_VHom}).

\begin{thm}\label{Thm_intoVhom}
Let $K$ be a triangulated homology $n$-manifold. Then, for every $\epsilon\in \mC^0(K)$, $V_\epsilon$ is a triangulated homology $(n-1)$-manifold. Moreover, if $K$ is a $PL$-triangulation of a topological $n$-manifold, then for every $\epsilon \in \mC^0(K)$, $V_\epsilon$ is a $PL$-triangulated topological $(n-1)$-manifold.
\end{thm}

 Recall that a homology $n$-manifold is a topological space $X$ such that for every point $x\in X$, the relative homology
 $H_*(X,X\setminus\{x\};\Z)$ is isomorphic to $H_*(\R^n,\R^n\setminus\{0\};\Z)$. Any smooth or topological $n$-dimensional manifold  is thus a homology $n$-manifold.

Poincar\'e duality holds true in such compact homology manifolds, see \cite{M}. 
And a triangulation $K$ is called piecewise linear ($PL$) if for every simplex $\tau\in K$, the link $\Lk(\tau,K)=\{\sigma\in K | \exists \eta \in K \mbox{ such that } \sigma,  \tau<\eta \mbox{ and }  \tau \cap \sigma=\emptyset\}$ is homeomorphic to a sphere, see Section~\ref{Sect_$Vepsilon$}. 

Our main goal is to understand the topology of $V_\epsilon$ when $\epsilon\in \mC^{k-1}(K)$ is chosen at random.
More precisely, for every $d\geq0$, let us denote by $\Sd^d(K)$ the $d^{th}$ barycentric subdivision of $K$, with the convention that $\Sd^0(K)=K$. 
When $K$ is finite, the asymptotic behavior of the number of simplices of $\Sd^d(K)$ in each dimension $p\in \{0,1,\ldots, n\}$ has been studied  in \cite{BW}, \cite{DPS}, see also \cite{SW1}. This number $f_p^d(K)=f_p(\Sd^d(K))$ is equivalent to $q_{p,n} f_n(K)(n+1)!^d$  as $d$ grows to $+\infty$ for some universal constant $q_{p,n}>0$, where $f_p(K)$ denotes the number of $p$-dimensional simplices of $K$. 
Let $\nu\in [0,1]$. For every $d\geq 0,$ we equip  $\mC^{k-1}(\Sd^d(K))$  with the product probability measure $\mu_\nu$ so that for every $\epsilon\in \mC^{k-1}(\Sd^d(K))$, the probability that  $\epsilon$ takes the value 0 on a $(k-1)$-simplex of $\Sd^d(K)$  is $\nu$, while the probability that it takes the value 1 is  $1-\nu$.  
When $K$ is finite, we set  $\E_{\nu,d}(\chi)=\int_{ \mC^{k-1}(\Sd^d(K))} \chi(V_\epsilon)d\mu_\nu(\epsilon)$ and for every $0\leq i\leq n-k,$
$\E_{\nu,d}(b_i)=\int_{ \mC^{k-1}(\Sd^d(K))} b_i(V_\epsilon)d\mu_\nu(\epsilon),$
  where $b_i(V_\epsilon)= \dim H_i(V_\epsilon)$ denotes the $i^{th}$ Betti number of $V_\epsilon$ with $\Z/2\Z$-coefficients and $\chi(V_\epsilon)=\sum_{i=0}^{n-k}(-1)^ib_i(V_\epsilon)$ its Euler characteristic. To simplify the notation, we will write $\E_{\nu, d}=\E_\nu$ since there won't be any ambiguity on the simplicial complex concerned. When $\nu=\frac{1}{2}$, we will moreover omit  $\nu$ from the notation.
Our main result is the following, see Corollaries~\ref{Cor_EbiEuler$k>1$} and  \ref{Cor_Inf}.

\begin{thm}\label{Thm_Intro} Let $k\in \{1,\ldots, n\}$. For every   $0\leq i\leq n-k$, there exist universal constants $c^+_{i}(n,k), c^-_{i}(n,k)>0, $ such that for every finite $n$-dimensional simplicial complex $K$ and every $0\leq i\leq n-k$,
$$c^-_{i}(n,k)\leq \liminf_{d\to +\infty}\frac{\E_\nu(b_i)}{(n+1)!^d f_n(K)}\leq \limsup_{d\to +\infty} \frac{\E_\nu(b_i)}{(n+1)!^d f_n(K)}\leq c^+_{i}(n,k).$$

Moreover,  $$\lim_{d\to +\infty}\frac{\E_\nu(\chi)}{(n+1)!^d f_n(K)}=\sum\limits_{i=0}^{n-k}(-1)^i c^+_{i}(n,k)$$ and 
$$\limsup_{d\to +\infty}\frac{1}{(n+1)!^df_n(K)}\big( \sum\limits_{j=0}^{i}(-1)^{i-j}\E_\nu(b_j)\big)\leq \sum\limits_{j=0}^{i}(-1)^{i-j} c^+_{j}(n,k).$$

\end{thm}

The universal constants $c^+_{i}(n,k)$ are linear combinations of the coefficients $q_{p,n}$ above, see   Definition~\ref{Defn_ $c_i^+(n,k)$}. 
The upper estimates in Theorem~\ref{Thm_Intro} actually follow from a computation of the expected face polynomial of the subcomplexes. The face polynomial  of a finite $n$-dimensional  simplicial complex $K$ is the polynomial $q_K(T)=\sum_{i=0}^{n}f_i(K)T^i.$ We denote by $\E_\nu(q_{V_\epsilon}(T))=\int_{C^{k-1}(K)}q_{V_\epsilon}(T)d\mu_\nu(\epsilon)$ the expected value of the face polynomial of a random subcomplex. Then, we prove the following (see Theorem~\ref{Thm_Enu}).

\begin{thm}\label{Intro_thmEnu}
Let $K$ be a finite $n$-dimensional simplicial complex and $k\in\{ 1,\ldots, n\}$. Then, for every $\nu\in[0,1]$, $$\E_{\nu}(q_{V_\epsilon}(T))=\int_Kq_{D(\sigma)}(T)dm_k(\sigma).$$

Moreover, if $K$ is a compact triangulated homology $n$-manifold and $k=1$, then $$\E_\nu(R_{V_\epsilon}(-1-T))=(-1)^{n} \E_\nu({R_{V_\epsilon}(T)}),$$ where $R_{V_\epsilon}(T)=Tq_{V_\epsilon}(T)-\chi(V_\epsilon)T.$
\end{thm}

Recall that the first barycentric subdivision $\Sd(K)$ inherits a decomposition into blocks, the block $D(\sigma)$ dual to $\sigma\in K$, see Section~\ref{Sect_$Vepsilon$} and \cite{M}. The right hand side in Theorem~\ref{Intro_thmEnu} is thus the total face polynomial of these blocks with respect to the measure $m_k$ on $K$. This measure $m_k$ equals to $\sum_{\sigma\in K}\mu_\nu\{\epsilon\in C^{k-1}(K)|\hat{\sigma}\in V_\epsilon\}\delta_{\hat{\sigma}}$, where $\hat{\sigma}$ is the barycenter of $\sigma$ and $\delta_{\hat{\sigma}}$ is the Dirac measure on it, see Section~\ref{Sect_mk} and \cite{SW1} for a study of such measures and integrals. This means that the density of $m_K$ at $\sigma\in K$ with respect to the canonical measure $\sum_{\sigma\in K}\delta_{\hat{\sigma}}$ studied in \cite{SW1} is given by the probability that $\hat{\sigma}$ belongs to $V_\epsilon$. When $K$ is a compact triangulated homology $n$-manifold, we obtain the analog of Theorem~\ref{Intro_thmEnu} where the face polynomial is replaced by the (simpler) block polynomial, see Theorem~\ref{Thm_vfi}.

 The constants  $c^-_{i}(n,k)$ depend on  some combinatorial complexity of the codimension $k$ closed submanifolds $\Sigma$ of $\R^n$, see Definition~\ref{Defn_$c_i^-$}. Namely, for every such closed connected codimension $k$ submanifold $\Sigma$ of $\R^n$, we define its complexity as the smallest value $m$ such that $(\R^n, \Sigma)$ gets homeomorphic to $(\stackrel{\circ}{\Delta}_n, {V}_\epsilon)$  for some $\epsilon \in \mC^{k-1}(\Sd^m(\Delta_n))$, where $\Delta_n$ denotes the standard simplex of dimension $n$. We actually estimate from below the asymptotic expected number of connected components  of $V_\epsilon$ which are homeomorphic to a given codimension $k$  closed connected submanifold $\Sigma$ of $\R^n$, see Theorem~\ref{Thm_cSigma}.  The complexity of surfaces in $\R^3$ is studied in Section~\ref{Sect_Complexity}, see Theorem~\ref{Thm_Complexity}.

When $k=1$, the subcomplexes $V_\epsilon$ turn out to inherit an additional $\CW$-complex structure, see Corollary~\ref{Cor_VCW}, which make it possible to improve the upper estimate in Theorem~\ref{Thm_Intro}. We prove the following, see  Corollaries~\ref{Cor_$E(bi)$} and \ref{Cor_E(chi)1}.
\begin{thm}\label{Thm_Intro$k=1$}
Let $K$ be a finite $n$-dimensional simplicial complex, $\nu\in [0,1]$, $k=1$ and $i\in \{0,\ldots, n-1\}$. Then,
$\E_\nu(b_i)\leq f_{i+1}(K)(1-\nu^{i+2}-(1-\nu)^{i+2}).$  Moreover,
$$\E_\nu(\chi)+\chi(K)=\sum_{i=0}^n(-1)^i(\nu^{i+1}+(1-\nu)^{i+1})f_i(K).$$
\end{thm}

Theorem~\ref{Thm_Intro$k=1$}
has the following quite surprising corollary which has already been observed by T.~Akita~\cite{A} with a different (non-probabilistic) method, but also follows from the symmetry property observed by I.G.~Macdonald~\cite{Mac}, see \cite{SW1}.

 \begin{cor} [\cite{A}] \label{Cor_Intro}If $K$ is a triangulated compact homology $2n$-manifold, then 

$$\chi(K)=\sum\limits_{p=0}^{2n}(-\frac{1}{2})^pf_p(K).$$
\end{cor}

We checked in \cite{SW1} that $-1$ together with $-\frac{1}{2}$ in even dimensions are the only universal parameters for which the polynomial $q_K(T)=\sum_{p=0}^{n}f_p(K)T^p$  equals $\chi(K)$ for every triangulated manifold $K$.
The paper is organized as follows. Section~\ref{Sect_$Vepsilon$} is devoted to a study of the topological structures of $V_\epsilon$. We prove in particular Theorem~\ref{Thm_intoVhom} and exhibit an additional $\CW$-structure on $V_\epsilon$ when $k=1$, see Corollary~\ref{Cor_VCW}.
   Section~\ref{Sect_mk} is devoted to a study of the measure $m_k$ and several computations, in particular  when $\nu=\frac{1}{2}$, see Corollary~\ref{Cor_MuZ}. In Section~\ref{Ssect_$k=1$}, we take profit of the $\CW$-complex structure to prove Theorem~\ref{Thm_Intro$k=1$} and Corollary~\ref{Cor_Intro}. Section~\ref{Ssect_$k>1$}  is devoted to the upper estimates in Theorem~\ref{Thm_Intro} in the general case ($k>1$), while Section~\ref{Ssect_HomMfd} is devoted to the special case of compact homology manifolds.
 
Finally,  we prove in Section~\ref{Sect_Lower} 
the lower estimates of Theorem~\ref{Thm_Intro}, see Theorem~\ref{Thm_cSigma}, Corollary~\ref{Cor_Inf}, and study  in the second part of Section~\ref{Sect_Lower} the combinatorial complexity of surfaces in $\R^3$, see Definition~\ref{Defn_Complexity}, Theorem~\ref{Thm_Complexity}.

These results thus provide counterparts in this combinatorial framework to the ones obtained in  \cite{GW15} and \cite{GW14, GW162} for the expected Betti numbers of real algebraic submanifolds of real projective manifolds or nodal domains in smooth manifolds respectively. The paper ends with an appendix devoted to a further  study and interpretation of the constants $c_i^+(n,k)$.

\vspace{0.5 cm}
\textbf{Acknowledgement :}
The second author is partially supported by the ANR project MICROLOCAL (ANR-15CE40-0007-01).

\section{Structure of the subcomplexes $V_\epsilon$}\label{Sect_$Vepsilon$}

The aim of this section is to prove that when $k=1$, the subcomplexes  $V_\epsilon$ inherit an additional 
$\CW$-complex structure. Moreover, they are homology manifold when $K$ itself is a homology manifold, see Corollary~\ref{Cor_VCW} and \ref{Cor_VHom}.

Let us start with recalling some definitions concerning simplicial complexes. 
Let $K$ be a finite simplicial complex and $\tau\in K$.
The \emph{star} of $\tau$ in $K$, denoted $\St(\tau, K)$, is the union of interiors of all simplices of $K$ having $\tau$ as a face. The closure $\overline{\St}(\tau, K)$ of $\St(\tau,K)$ is the union of all simplices of $K$ having $\tau$ as a face.
The \emph{link} of $\tau$ in $K$, denoted $\Lk(\tau,K)$, is the union of all simplices of $K$ lying in $\overline{\St}(\tau, K)$ that are disjoint from $\tau$. 

Let us recall as well that the \emph{join} $K\ast L$ of the simplicial complexes $K$ and $L$ is the simplicial complex 
whose simplices are the joins $\sigma_K\ast\sigma_L$ where $\sigma_K$ (respectively $\sigma_L$) are the simplices of $K$ (respectively $L$), including $\emptyset$.
If $\sigma_K=[e_0,\ldots, e_k]$ and $\sigma_L=[f_0,\ldots,f_l]$ then by definition $\sigma_K\ast \sigma_L=[e_0,\ldots,e_k,f_0,\ldots,f_l]$. In particular, if $L$ is a point, then $K\ast L$ is the cone over $K$ centered at $L$, see~\cite{M}.

\begin{thm}\label{Thm_Link}
Let $K$ be a locally finite $n$-dimensional simplicial complex. For every $1\leq k\leq n$, every $\epsilon \in \mC^{k-1}(K)$ and every $\sigma=[\hat{\sigma}_0,\ldots, \hat{\sigma}_p]\in V_\epsilon\subset \Sd(K)$,
there exists a canonical isomorphism of simplicial complexes between the link  $\Lk(\sigma, V_\epsilon)$ of $\sigma$ in $V_{\epsilon}$ and $(V_\epsilon\cap \Bd \sigma_0)\ast\Sd(\Bd\Lk(\sigma_0,\sigma_1))\ast\dots\ast\Sd(\Bd\Lk(\sigma_{p-1},\sigma_{p}))\ast\Sd(\Lk(\sigma_p,K)).$
\end{thm}

\begin{pf}
By definition, a simplex $\tau$ of $\overline{\St}(\sigma, V_\epsilon)$ is of the form $[\hat{\tau_0},\ldots,\hat{\tau}_q]$  such that there exists a sequence of subindices $0\leq i_0<i_1<\dots< i_p\leq q$ satisfying $\tau_{i_j}=\sigma_j$, where $\tau_0<\tau_1<\dots<\tau_q$  and $\epsilon$ restricted to ${\tau_0}$ is not constant. Therefore,  a simplex of $\Lk(\sigma, V_\epsilon)$ is of the form   $[\hat{\tau_0},\ldots,\hat{\tau}_{i_{0}-1}, \hat{\tau}_{i_{0}+1}\ldots\hat{\tau}_{i_{p}-1}, \hat{\tau}_{i_{p}+1}\ldots\hat{\tau}_q]$, and thus it has a canonical decomposition $[\hat{\tau}_0,\ldots,\hat{\tau}_{{i_0}-1}]\ast [\hat{\tau}_{{i_0}+1},\ldots,\hat{\tau}_{{i_1}-1}]\ast\dots\ast[\hat{\tau}_{{i_p}+1},\ldots,\hat{\tau}_q]$. 
Note that the first term $[\hat{\tau}_0,\ldots,\hat{\tau}_{{i_0}-1}]$ is an element of $V_\epsilon\cap \Bd\sigma_0$ as  $\tau_{i_0-1}$ is a proper face of $\sigma_0$ and $\epsilon$ being non-constant on $\tau_0$ is the condition of belonging to $V_\epsilon$. The last term $[\hat{\tau}_{{i_p}+1},\ldots,\hat{\tau}_q]$ can be identified to $[\widehat{\Lk(\sigma_p,\tau_{{i_{p}}+1})},\ldots,\widehat{\Lk(\sigma_p,\tau_q})]$ by replacing each $\hat{\tau}_s$ with $\widehat{\Lk(\sigma_p,\tau_s)}$ for every $s\in\{{i_p}+1,\ldots, q\}$, while for every $j\in \{0,\ldots,{p-1}\}$, the intermediate term $[\hat{\tau}_{{i_j}+1},\ldots,\hat{\tau}_{{i_{j+1}}-1}]$ can be identified to $[\widehat{\Lk(\sigma_{j},\tau_{{i_j}+1})},\ldots,\widehat{\Lk(\sigma_{j},\tau_{{i_{j+1}}-1}})]$, respectively. Therefore, doing so, we obtain an element of $\Sd(\Lk(\sigma_p,K))$ and respectively an element of $\Sd(\Bd\Lk(\sigma_j,\sigma_{{{j+1}}}))$ for every $j\in \{0,\ldots,p-1\}$.

Let $L$ denote the join $(V_\epsilon\cap \Bd\sigma_0) \ast\Sd(\Bd\Lk(\sigma_0,\sigma_1))\ast\dots\ast\Sd(\Bd\Lk(\sigma_{p-1},\sigma_{p}))\ast\Sd(\Lk(\sigma_p,K))$,  we  then obtain a canonical simplicial map 

$$\begin{array}{rcl}
\Phi :\Lk(\sigma, V_\epsilon)&\to&L\\

[\hat{\tau}_0,\ldots,\hat{\tau}_q]&\to& [\hat{\tau}_0,\ldots,\hat{\tau}_{{i_0}-1}]\ast [\widehat{\Lk(\sigma_0,\tau_{{i_{0}}+1})},\ldots,\widehat{\Lk(\sigma_0,{\tau_{{i_1}-1}}})]\ast\ldots\ast\\
&&\ldots\ast[\widehat{\Lk(\sigma_{p},\tau_{{i_{p}}+1})},\ldots,\widehat{\Lk(\sigma_p,\tau_q})].
\end{array}$$

To define the inverse map, we note that  an element of $V_\epsilon\cap \Bd\sigma_0$ is of the form $\Gamma=[\hat{\gamma}_0,\ldots,\hat{\gamma}_{{i_0}-1}]$ such that $\gamma_0<\dots<\gamma_{{i_0}-1}<\sigma_0$ and that $\epsilon|_{\gamma_0}$ is not constant. 
Meanwhile, for every  $j\in \{0,\ldots,{{p-1}}\}$, an element of $\Sd(\Bd\Lk(\sigma_j,\sigma_{j+1}))$ is of the form 
$\Gamma_{j}=[\hat{\gamma}_{i_j+1},\ldots,\hat{\gamma}_{{i_{j+1}}-1}]$ such that $\gamma_{i_{j}+1}<\dots<\gamma_{i_{j+1}-1}<\sigma_{j+1}$ and that $\gamma_s$ is disjoint from $\sigma_{j}$ for every $s\in\{{i_j}+1,\ldots, i_{j+1}-1\}$.
Finally,  an element of  $\Sd(\Lk(\sigma_p,K))$ is of the form $\Gamma_{{p}}=[\hat{\gamma}_ {i_{p}+1},\ldots,\hat{\gamma}_q]$ such that $\gamma_ {i_{p}+1}<\dots<\gamma_q\in \St(\sigma_p, K)$ and that $\gamma_s$ is disjoint from $\sigma_p$ for every $s\in\{{i_{p}+1},\ldots,q\}$.
Therefore, there is a canonical simplicial map $\Psi :L\to\Lk(\sigma, V_\epsilon)$ such that 
$$\begin{array}{rcl}
\Psi :L&\to&\Lk(\sigma, V_\epsilon)\\
\Gamma\ast\Gamma_0\ast\ldots\ast\Gamma_p

&\to&
[\hat{\gamma}_0,\ldots,\hat{\gamma}_{{i_0}-1}, \widehat{\sigma_0\ast\gamma_{{i_0}+1}},\ldots,\widehat{\sigma_0\ast\gamma_{{i_1}-1}}, \ldots,\\
&&\\
&&\ldots,\widehat{\sigma_{p-1}\ast\gamma_{{i_{p}}-1}}, \widehat{\sigma_p\ast\gamma_{{i_p}+1}},\ldots,\widehat{\sigma_p\ast\gamma_q}]. 
\end{array}$$

The maps $\Phi$ and $\Psi$ are inverse to each other, hence the result.
\end{pf}

The case  $k=1$ is of special interest due to the following proposition and its corollaries. 
\begin{prop}\label{Prop_$k=1$}
Let $K$ be a locally finite simplicial complex. For every $\epsilon \in \mC^0(K)$, 
the intersection of $V_\epsilon$ with a $p$-simplex $\sigma$ of $K$, if not empty, is isotopic to an affine hyperplane section $H_\sigma\subset \sigma$ which separates the vertices of  $\sigma$ labelled 1 from the vertices labelled 0. 
\end{prop}
By isotopic in Proposition~\ref{Prop_$k=1$} we mean that there exists a continuous family of homeomorphisms of $\sigma$ from the identity to a homeomorphism which maps $V_\epsilon\cap \sigma$ to $H_\sigma$.

\begin{figure}[h]
   \begin{center}
    \includegraphics[scale=0.4]{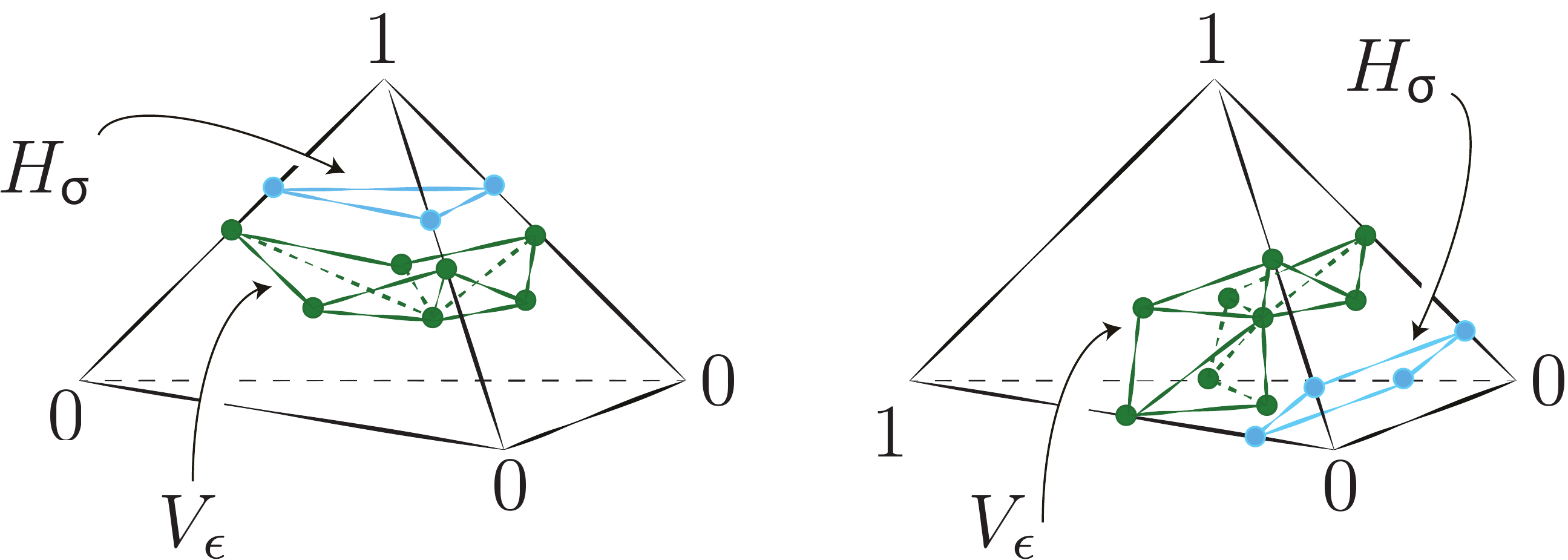}
    \caption{Examples of $V_\epsilon\cap \sigma$ versus $H_\sigma$ for $p=3$.}
    \label{VepsilonvsH}
      \end{center}
 \end{figure}

Figure~\ref{VepsilonvsH} exhibits  some examples of hyperplane sections given by Proposition~\ref{Prop_$k=1$}.

\begin{rem}\label{Rem_Ex}
Note that  Proposition~\ref{Prop_$k=1$} does not hold true for $k>1$. Indeed, consider for example $\epsilon \in \mC^{1}(T)$ where $T$ is a tetrahedron such that $\epsilon$ has constant value 1 on all the edges of $K$, the case depicted on the right of Figure~\ref{Vepsilon}.
Then, $V_\epsilon$ is the cone over the barycenters of  the four 2-dimensional faces of $T$ centered at the barycenter of $T$. \end{rem}

\begin{pf}  Suppose that  $\epsilon$ takes the value 1 on $0<j+1<p+1$  vertices of a $p$-simplex $\sigma$. 
There is then a $j$-face $\tau_1$ of $\sigma$ with all its vertices labelled 1 and a $(p-j-1)$-face $\tau_0$ with vertices labelled by 0. 
Let $H_\sigma\subset \sigma$ be a hyperplane section separating $\tau_0$ from $\tau_1$ intersecting transversally exactly $(j+1)(p-j)$-many 1-faces which are neither in $\tau_0$ nor in $\tau_1$.
 
Let us prove by induction on the dimension of the skeleton of  $\sigma$ that $V_\epsilon\cap \sigma$ is isotopic to $H_\sigma$.
As $\epsilon$ takes value $1$ on exactly $j+1$ vertices, there are $(j+1)(p-j)$-many edges on which $d\epsilon$ is non trivial. Those are the edges which are neither on $\tau_0$ nor on $\tau_1$.
Thus the intersection of $V_\epsilon$ with the 1-skeleton of $\sigma$ is a set of isolated points, one point on each of these $(j+1)(p-j)$-many edges.
So, on each such edge we have two points, one defined by the intersection of $V_\epsilon$ and the other by the intersection of $H_\sigma$.
We perform an isotopy which takes one set of points to the other.
Now, let us suppose that $V_\epsilon$ can be isotoped to $H_\sigma$ on the $l$-skeleton for $l\geq 1$ and choose such an isotopy. Let $\eta$ be a $(l+1)$-face of $\sigma$.
The intersection $V_\epsilon \cap \eta$, if not empty, is the cone centered at the barycenter of $\eta$ over the intersection of $V_\epsilon$ with the boundary of $\eta$. 
Besides $H_\sigma \cap \eta$ is an affine hyperplane section. Let $b$ be a point in the hyperplane section so that $H_\sigma \cap \eta$ is the cone over $H_\sigma \cap \partial \eta$ centered at $b$. 
By means of choosing a path $(b_t)_{t\in[0,1]}$ in $\stackrel{\circ}{\sigma}$ from $b$ to the barycenter of $\eta$, we can get an isotopy between $H_\sigma\cap \eta$ and $V_\epsilon \cap \eta$ by taking the cone centered at $b_t$ over the isotopy between $H_\sigma\cap \partial \eta$ and $V_\epsilon\cap \partial \eta$.
\end{pf}

\begin{rem} \label{Rem_Prop_$k=1$}
 Note that when $K$ is a convex triangulation of the moment polytope of some toric manifold, then $\epsilon \in \mC^0(K)$ 
provides a distribution of signs  on every vertex of $K$ and the collection $(H_\sigma)_{\sigma_\in K}$ defines  on $K$ the (tropical) hypersurface constructed by O.~Viro~\cite{V1} in his patchwork theorem, provided $H_\sigma$ is chosen to intersect every edge in its middle point.  From Proposition~\ref{Prop_$k=1$}, we thus deduce that the pair $(K, V_\epsilon)$  is homeomorphic  to the pair $(K, H)$ defined by O.~Viro, where $H=\cup_{\sigma\in K} H_\sigma$. 
\end{rem}

\begin{cor}\label{Cor_VCW}
Let $K$ be a locally finite $n$-dimensional simplicial complex and $\epsilon\in \mC^0(K)$. Then, $V_\epsilon$ inherits the structure of a  $\CW$-complex, having a cell of dimension $p$ for every $(p+1)$-simplex of $K$ on which $\epsilon$ is not constant, $p\in\{0,\ldots, n-1\}$.
\end{cor}
\begin{pf} From Proposition~\ref{Prop_$k=1$}, we indeed know that
the intersection of $V_\epsilon$ with any $(p+1)$-simplex of $K$ on which $\epsilon$ is not constant is homeomorphic to a $p$-cell. Hence, the result.
\end{pf}

\begin{cor} \label{Cor_VHom}Let $K$ be a  triangulated homology $n$-manifold. Then, for every $\epsilon \in \mC^0(K)$, $V_\epsilon$  is a triangulated $(n-1)$-homology manifold. Moreover, if $K$ is a $PL$-triangulation of a homology $n$-manifold, then for every $\epsilon\in \mC^0(K)$, $V_\epsilon$ is a $PL$-triangulation of a homology $(n-1)$-manifold. 
\end{cor}

\begin{pf} From Theorem~\ref{Thm_Link} it follows that for every $\sigma\in V_\epsilon$, the link $\Lk(\sigma, V_\epsilon)$ is canonically isomorphic  to $(V_\epsilon\cap \Bd\sigma_0)\ast\Sd(\Bd\Lk(\sigma_0,\sigma_1))\ast\ldots\ast\Sd(\Bd\Lk(\sigma_{p-1},\sigma_{p}))\ast\Sd(\Lk(\sigma_p,K))$, where the first term is homeomorphic to a sphere by Proposition~\ref{Prop_$k=1$} and the intermediate terms are by definition homeomorphic to spheres.
Finally, the last term is a homology sphere by Lemma~63.1 of \cite{M} in the case where $K$ is a triangulated homology manifold (respectively, it is a sphere in the case where $K$ is a $PL$-triangulation). Therefore, the link of any simplex of $V_\epsilon$ is the join of a homology sphere with spheres, which is a homology sphere (respectively, join of spheres which is a sphere).  Hence the result.
\end{pf}

\begin{ex} \label{Ex_Notmfd}If $K$ is a triangulation of a closed manifold, then $V_\epsilon$ may not be a triangulation of a submanifold.
A counterexample can be constructed from the double suspension of the Poincar\'e sphere. Namely, let $K$ be a triangulation of the Poincar\'e sphere and let $S$ denote the simplicial complex of a 0-dimensional sphere. We take the double suspension ($S\ast S\ast K$) of $K$ together with a simplicial complex structure obtained by considering successive cones first over the simplexes of $K$ centered at vertices of $S$ and then over $S\ast K$ centered at the vertices of $S$. The obtained complex, denoted $\widetilde{K}$, is a triangulation of a 5-dimensional sphere, see \cite{Can, Ed}. Now, let us consider $\epsilon\in \mC^0(\widetilde{K})$ such that $\epsilon$ takes value 1 on one of the four vertices that corresponds to one of the four centers of suspensions and zero on all other vertices. There is a natural isotopy from
$V_\epsilon$ to $S\ast K$. The latter is not a submanifold as the link of the two points corresponding to the center of suspension are Poincar\'e spheres.
\end{ex}

\begin{rem} \label{Rem_VHom}
Corollary~\ref{Cor_VHom} does not hold true for $k>1$. Indeed, in the case of the tetrahedron discussed in Remark~\ref{Rem_Ex} for example, the link $\Lk(\widehat{T}, V_\epsilon)$ of the barycenter $\widehat{T}$ of $T$ is the set of four vertices which is not a homology sphere.  This example can be implemented in any triangulated homology $3$-manifold $K$ so that   $V_\epsilon$ need not be a triangulated homology manifold, although $K$ is.
\end{rem}

\begin{cor} Let $K$ be an even dimensional compact triangulated homology manifold. Then, for every $\epsilon\in  \mC^0(K), \,\chi(V_\epsilon)=0$. 
\end{cor}

\begin{pf} By Corollary~\ref{Cor_VHom}, $V_\epsilon$ is a homology manifold in which case the Poincar\'e duality  with $\Z/2\Z$-coefficients applies, see Chapter~8 of \cite{M}. When the dimension of $K$ is even, the dimension of $V_\epsilon$ is odd, hence the result.
\end{pf}

\section{The induced measures $m_k$}\label{Sect_mk}
Let $K$ be a locally finite $n$-dimensional simplicial complex and $k\in\{1,\ldots,n\}$, for every $\sigma\in K$, we set

$$m_k(\sigma)=\mu_\nu\{\epsilon\in C^{k-1}(K, \Z/2\Z) | \hat{\sigma}\in V_{\epsilon}\}.$$
This is the probability that the barycenter $\hat{\sigma}$ belongs to $V_\epsilon$. It defines a measure on $K$, namely

$$m_k=\sum_{\sigma\in K}\mu_\nu\{\epsilon \in C^{k-1}(K,\Z/2\Z)| \hat{\sigma}\in V_\epsilon\}\delta_{\hat{\sigma}},$$
where $\delta_{\hat{\sigma}}$ denotes the Dirac measure on $\hat{\sigma}$. We likewise set, for every $p\in\{0,\ldots, n\},$
$$\gamma_{p,K}=\sum_{\sigma\in K^{[p]}}\delta_{\hat{\sigma}},$$ where $K^{[p]}$ denotes the set of $p$-dimensional simplices of $K$.
 
We proved in  \cite{SW1} that $\frac{1}{(n+1)!^d}\gamma_{p,\Sd^d(K)}$ weakly converges to $q_{p,n}\dvol_{K}$ as $d$ grows to $+\infty$, where $q_{p,n}>0$ is some universal constant and $\dvol_{K}=\sum_{\sigma\in K^{[n]}}(f_\sigma)_*(\dvol_{\Delta_n})$.
In the latter $\dvol_{\Delta_n}$ is the  Lebesgue measure of the standard simplex $\Delta_n$ normalized in such a way that it has total measure 1 and $f_\sigma:\Delta_n\to \sigma$ is some affine isomorphism. Our aim in this section is to study the measure $m_k$.
\begin{prop} \label{Prop_mk} Let $K$ be a locally finite $n$-dimensional simplicial complex and $k\in\{1,\ldots, n\}.$ Then, $m_k=\sum_{p=k}^{n}\big(1-\mu_\nu(Z^{k-1}(\Delta_p))\big)\gamma_{p,K}$, where $Z^{k-1}(\Delta_p)\subset C^{k-1}(\Delta_p)$ denotes the subspace of $(k-1)$-cocycles of $\Delta_p$ with $\Z/2\Z$-coefficients.
\end{prop}    

\begin{pf} By definition, for every $\sigma\in K$ and every $\epsilon \in C^{k-1}(K)$, $\hat{\sigma}$ belongs to $V_\epsilon$ if and only if the restriction of $d\epsilon$ to $\sigma$ does not vanish, that is $\epsilon|_\sigma\notin Z^{k-1}(\sigma).$ Since $\mu_\nu$ is a product measure, we deduce that $\mu_\nu\{\epsilon \in C^{k-1}(K)|\hat{\sigma}\in V_\epsilon\}=1-\mu_\nu(Z^{k-1}(\sigma)).$ Moreover, by definition of $\mu_\nu$, $\mu_\nu(Z^{k-1}(\sigma))$ only depends on the dimensions of $\sigma$. Finally, $m_k=\sum_{\sigma\in K}\big(1-\mu_\nu(Z^{k-1}(\sigma))\big)\delta_{\hat{\sigma}}=\sum_{p=k}^{n}\big(1-\mu_\nu(Z^{k-1}(\Delta_p))\big)\gamma_{p,K}.$
\end{pf}

We may  additionally set $m_0=\sum_{p=0}^{n}\gamma_{p,K}$ so that all statements involving $m_k$ will make sense when $k=0$ as well, but this case is of no interest for us.

We have been able to compute explicitly  the universal constants $\mu_\nu(Z^{k-1}(\Delta_p))$ appearing in Proposition~\ref{Prop_mk}  in several cases, in particular, when  $\nu=\frac{1}{2}$, see Corollary~\ref{Cor_MuZ}. These computations are based on the following theorem.

\begin{thm}\label{Thm_MuZ}
For every $1\leq k\leq p$,
$$
\mu_\nu(Z^{k-1}(\Delta_p))=\begin{cases} 
\nu^{p+1}+(1-\nu)^{p+1} & \mbox{if } k=1,\\
 \int_{\eta \in \mC^{k-2}(\Delta_{p-1})}\mu_\nu(d\eta)d\mu_\nu(\eta)& \mbox{if } k\geq 2.\\
  \end{cases}
$$
\end{thm}

The case of $k\geq2$ in Theorem~\ref{Thm_MuZ} is thus the expected value of the random variable $\mu_\nu\circ d$ on the probability space $C^{k-2}(\Delta_{p-1}).$

\begin{pf}
When $k=1$, $Z^0(\Delta_p)=\{0,1\}$ where 0 (respectively 1)
denotes the 0-cochain which is constant and equal to 0 (respectively 1) on $C_0(\Delta_p)$. By definition of $\mu_\nu$, $\mu_\nu(0)=\nu^{p+1}$ and $\mu_\nu(1)=(1-\nu)^{p+1}$ since $\Delta_p$ has $p+1$ vertices. Let us assume now that $k\geq 2$ and let $s$ be a vertex of $\Delta_p$. We are going to prove that whatever the value of $\epsilon$ on the $(k-1)$-faces of $\Delta_p$ containing $s$ is, there is a unique way to extend $\epsilon$  to all $(k-1)$-faces of $\Delta_p$ in such a way that $d\epsilon$ restricted to $\Delta_p$ vanishes. Indeed, let us assume that we fix a value of $\epsilon$ on the $(k-1)$-faces of $\Delta_p$ containing $s$ and consider a $k$-face $\tau$ of $\Delta_p$ that contains $s$. In this case, among the $k+1$ many $(k-1)$-faces of $\tau$ there is only one, say $\eta$, which does not contain $s$.  By the assumption above, all $(k-1)$-faces of $\tau$ but $\eta$ are labelled by $\epsilon$. There is a bijection on those labelled $(k-1)$-faces of $\tau$ and the  $(k-2)$-faces of $\eta$, since the former are  cones over the latter.
This bijection thus induces labels on the $(k-2)$-faces of $\eta$. Let us assign to $\eta$ the value 1 if an odd number of its codimension-1 faces are labelled 1; and 0 otherwise.
By doing so, in either case an even number of $(k-1)$-faces of $\tau$  get labelled 1, which results in $<d\epsilon,\tau>=0$. 
Moreover, this way  is the only way to label $\eta$ for having $<d\epsilon,\tau>=0$. 
At this point, $\epsilon$ has been extended to all $(k-1)$-faces of $\Delta_p$ and we have to check that $d\epsilon$ vanishes on $\Delta_p$.

Let then now $\tau$ be a $k$-face of $\Delta_p$ which does not contain $s$. By definition, the restriction of $\epsilon$ to $\tau$ equals $d\tilde{\epsilon}$ where $\tilde{\epsilon}\in \mC^{k-2}(\tau)$ is inherited by the values of $\epsilon$
on cones over the $(k-2)$-faces of $\tau$, centered at $s$.
Since $d^2=0$, we deduce that  $<d\epsilon, \tau>=<d^2\tilde{\epsilon}, \tau>=0$. 
Now,  we deduce that 

$$\mu_\nu(Z^{k-1}(\Delta_p))=\sum_{\eta \in \mC^{k-2}(\Lk(s,\Delta_p))} \mu_\nu(\eta) \mu_\nu(d \eta),$$
since $\mu_\nu$ is a product measure. Recall that
$\Lk(s,\Delta_p)$ denotes the link of $s$ in $\Delta_p,$ that is the $(p-1)$-simplex spanned by  all the vertices of $\Delta_p$ but $s$. 
 The coefficient $\mu_\nu(\eta)$ computes the value of $\mu_\nu$ on cones centered at $s$ over the $(k-2)$-faces of $\Lk(s, \Delta_p)$, while the coefficients $\mu_\nu(d\eta)$ computes the value of $\mu_\nu$ on the $(k-1)$-faces of $\Lk(s,\Delta_p).$ 
  \end{pf}

\begin{cor}\label{Cor_MuZ}
\begin{enumerate}

\item  If $\nu=\frac{1}{2}$, then for every $1\leq k\leq p,$ $$\mu_\nu(Z^{k-1}(\Delta_p))=\frac{1}{2^{\binom{p}{k}}}.$$

\item For every $\nu\in [0,1],$

$$\mu_\nu(Z^{k-1}(\Delta_p))=\begin{cases}
\frac{1}{2}(1+(2\nu-1)^{k+1}) & \mbox{if } 1\leq k=p,\\
&\\
\nu^{\binom{p+1}{2}}\sum\limits_{l=0}^{p}\binom{p}{l}\big(\frac{1-\nu}{\nu}\big)^{(l+1)(p-l)} & \mbox{if } 2=k\leq p.\\
\end{cases}$$
The first line remains valid for every $1\leq k\leq p$ if $\nu\in \{0,1\}$.
\end{enumerate}
\end{cor}

\begin{pf}
Let $\nu=\frac{1}{2}$. If $k=1$, the result directly follows from Theorem~\ref{Thm_MuZ}.  If $2\leq k\leq p$, note that the $(p-1)$-simplex has $\binom{p}{k-1}$ faces of dimension $k-2$ and   $\binom{p}{k}$ faces of dimension $k-1$, so we deduce from Theorem~\ref{Thm_MuZ} that 
 
   $$\begin{array}{lcl}
\mu_\nu(Z^{k-1}(\Delta_p))&=&\sum\limits_{\eta\in \mC^{k-2}(\Delta_{p-1})}\mu_\nu(\eta)\mu_\nu(d\eta)\\
&=&\sum\limits_{\eta\in \mC^{k-2}(\Delta_{p-1})}\frac{1}{2}^{\binom{p}{k-1}+\binom{p}{k}}\\
&=&\frac{1}{2^{\binom{p}{k}}}.
 \end{array}$$

 Now, let $\nu\in[0,1]$. If $k=1=p$, the result directly follows  from Theorem~\ref{Thm_MuZ}. If $2\leq k=p$, then from Theorem~\ref{Thm_MuZ},
 
 $$\begin{array}{lll}
\mu_\nu(Z^{k-1}(\Delta_p))&=&\sum_{\eta\in \mC^{k-2}(\Delta_{p-1})}\mu_\nu(\eta)\mu_\nu(d\eta)\\
 &&\\
 &=&\sum\limits_{l=0}^{\lfloor \frac{p}{2}\rfloor}\binom{p}{2l}(1-\nu)^{2l}\nu^{p-2l+1}+\sum\limits_{l=0}^{\lfloor \frac{p-1}{2}\rfloor}\binom{p}{2l+1}\nu^{p-2l-1}(1-\nu)^{2l+2}
 \end{array}$$
 since the $(p-1)$-simplex $\Delta_{p-1}$ has $p$ facets and the value of $d\epsilon$ on its unique $(p-1)$-face depends on the parity of the number of facets where $\epsilon=1.$
We deduce 
 $$\begin{array}{lll}
\mu_\nu(Z^{k-1}(\Delta_p))&=&\nu\sum\limits_{l=0}^{p}\binom{p}{l}(1-\nu)^{l}\nu^{p-l}(-1)^l+\sum\limits_{l=0}^{p}\binom{p}{l}(1-\nu)^{l}\nu^{p-l}(\frac{1-(-1)^l}{2})\\
 &=& \nu(2\nu-1)^p+\frac{1}{2}(1-(2\nu-1)^p)\\
 &=& \frac{1}{2}(1+(2\nu-1)^{p+1}).
  \end{array}$$
 
 When $\nu=1$, for every $2\leq k\leq p$, 
 $$\mu_\nu(Z^{k-1}(\Delta_p))=\mu_\nu(0)\mu_\nu(d0)=1,$$ since $\mu_\nu$ is then a Dirac on $0\in \mC^{k-2}(\Delta_{p-1}),$ while when $\nu=0,$
  $$\mu_\nu(Z^{k-1}(\Delta_p))=\mu_\nu(1)\mu_\nu(d1)=\begin{cases}0& \mbox{if $k$ is even, }\\ 1& \mbox{if $k$ is odd.}\\ \end{cases} $$
  The latter remains valid for $k=1$, from Theorem~\ref{Thm_MuZ}.
  
  Now if  $2=k\leq p$, then  from Theorem~\ref{Thm_MuZ},

  $$\begin{array}{lcl}
\mu_\nu(Z^{k-1}(\Delta_p))&=&\sum\limits_{\eta\in \mC^{0}(\Delta_{p-1})}\mu_\nu(\eta)\mu_\nu(d\eta)\\
&&\\
 &=&\sum\limits_{l=0}^p\binom{p}{l}\nu^l(1-\nu)^{p-l}\nu^{(\binom{l}{2}+\binom{p-l}{2})}(1-\nu)^{l(p-l)},
 \end{array}$$
 where $\binom{l}{2}=0$ (respectively $\binom{p-l}{2}=0$) if $l<2$ (respectively $p-l<2$). Indeed, the standard $(p-1)$-simplex has $p$ vertices and if $\epsilon$ takes the value 0 on $l$ of them, it has $\binom{l}{2}+\binom{p-l}{2}$ edges where $d\epsilon =0$ and $l(p-l)$ edges where $d\epsilon=1$ . The result follows from the relation $\binom{l}{2}+\binom{p-l}{2}=\binom{p}{2}-l(p-l)$, valid for every $l\in \{0,\ldots, p\}$.
\end{pf}

\section{Asymptotic topology of $V_\epsilon$, upper estimates}

\subsection{The case $k=1$}\label{Ssect_$k=1$}
Let $K$ be a finite $n$-dimensional simplicial complex. 
When $k=1$, $V_\epsilon$ admits a $\CW$-complex structure given by Corollary~\ref{Cor_VCW} whatever $\epsilon \in C^{k-1}(K)$ is.
Let $\tilde{f}_p(V_\epsilon)$ be the number of  $p$-cells of $V_\epsilon$, $p\in\{0,\ldots, n\}$, for this $\CW$ structure and ${\tilde{q}}_{V_\epsilon}(T)=\sum_{p=0}^n \tilde{f}_p (V_\epsilon) T^p$.

The expected value of this polynomial is given by the following theorem.

\begin{thm} \label{Thm_E(tilde{q})} Let ${K}$ be a finite $n$-dimensional simplicial complex and $\nu \in [0,1]$. Then, when $k=1$,
$$\E_{\nu}(T{\tilde{q}}_{\tiny{V_\epsilon}}(T))=\int_KT^{\dim\sigma}dm_1(\sigma)=q_{\tiny{K}}(T)-\nu q_K(\nu T)-(1-\nu)q_K((1-\nu)T),$$ where for every $\epsilon \in   \mC^0(K)$, $V_\epsilon$ is equipped with the $\CW$-complex structure given by Corollary~\ref{Cor_VCW}. In particular, for every $p\in \{0,\ldots,n-1\}.$

$$\E_\nu(\tilde{f}_p(V_\epsilon))=m_1(K^{[p+1]})=f_{p+1}(K)(1-\nu^{p+2}-(1-\nu)^{p+2}).$$
\end{thm}

\begin{pf} For every $\nu\in [0,1]$, $p\in \{0,\ldots,n-1\},$

$$\begin{array}{rcl} 
\E_{\nu}(\tilde{f}_p(V_\epsilon))&=&\int_{ \mC^0(K)}\tilde{f}_p(V_\epsilon) d\mu_\nu(\epsilon)\\
&&\\
&=&\sum\limits_{\epsilon \in  \mC^0(K)}  \#\{\sigma \in K^{[p+1]} |V_\epsilon \cap  \stackrel{\circ}{\sigma}\neq \emptyset \}\mu_\nu(\epsilon)\\
&&\\
&=&\sum\limits_{\sigma\in K^{[p+1]}}  \mu_\nu\{ \epsilon \in  \mC^0(K) | V_\epsilon \cap \stackrel{\circ}{\sigma} \neq \emptyset\}\\
&&\\
&=&m_1(K^{[p+1]})\\
&&\\
&=& f_{p+1}(K)(1-\nu^{p+2}-(1-\nu)^{p+2}),
\end{array}$$
since  $\mu_\nu$ is a product measure and the probability that $\epsilon$ is identically 0 (respectively 1) on the $p+2$ vertices of $\sigma$ is $\nu^{p+2}$ (respectively $(1-\nu)^{p+2})$.
This proves the second part  together with the first equality of the theorem.  Then,

$$\begin{array}{rcl} 
\E_{\nu}(T\tilde{q}_{{V_\epsilon}}(T))&=&\sum\limits_{p=0}^{n-1}f_{p+1}(K)(1-\nu^{p+2}-(1-\nu)^{p+2})T^{p+1}\\
&&\\
&=&\sum\limits_{p=0}^{n}f_{p}(K)(1-\nu^{p+1}-(1-\nu)^{p+1})T^{p}\\
&&\\
&=&q_K(T)-\nu q_K(\nu T)-(1-\nu)q_K((1-\nu)T).
\end{array}$$
\end{pf}

\begin{cor} \label{Cor_$E(bi)$}Let $k=1$. For every finite  $n$-dimensional simplicial complex $K$, every $\nu\in [0,1]$ and every $0\leq i \leq n-1$, $$\E_{\nu}(b_i)\leq m_1(K^{[i+1]})= f_{i+1}(K)(1-\nu^{i+2}-(1-\nu)^{i+2})$$ and  $$ \,\, \sum_{j=0}^{i} (-1)^{i-j} \E_{\nu}(b_j)\leq \sum_{j=0}^{i} (-1)^{i-j} {f}_{j+1}(K)(1-\nu^{j+2}-(1-\nu)^{j+2}),$$ where the latter is an equality if $i=n-1$.
\end{cor}

\begin{pf} By Corollary~\ref{Cor_VCW}, $V_\epsilon$ inherits a $\CW$-complex structure for all $\epsilon \in  \mC^0(K)$. From cellular homology theory, we deduce that $b_i(V_\epsilon)\leq \tilde{f}_i(V_\epsilon)$ for all $0\leq i \leq n-1$ and the Morse inequalities $\sum_{j=0}^{i} (-1)^{i-j} b_j(V_\epsilon)\leq \sum_{j=0}^{i} (-1)^{i-j} \tilde{f}_{j}(V_\epsilon)$.
The result follows by integrating over $ \mC^0(K)$ and applying the second part of Theorem~\ref{Thm_E(tilde{q})}.
\end{pf}

\begin{cor} \label{Cor_E(chi)1}Let $k=1$. For every finite simplicial complex $K$  and every  $\nu\in [0,1]$,
$$\E_{\nu}(\chi)+\chi(K)=\nu q_K(-\nu)+(1-\nu)q_K(\nu-1).$$ 
In particular, $\E(\chi)+\chi(K)= q_K(-\frac{1}{2})$.
\end{cor}

\begin{pf} The first part of the result follows from Theorem~\ref{Thm_E(tilde{q})}, by letting $T=-1$ and the second part by letting further $\nu=\frac{1}{2}.$
\end{pf}

\begin{cor} \label{Cor_R_nu} Let $k=1$ and $K$ be a triangulation of a compact homology $n$-manifold. For every $\nu\in [0,1],$
$$\E_{\nu}(\chi)=((-1)^n-1)R_K(-\nu).$$

In particular, if $n$ is odd, the values of the polynomial $q_n^\infty(T)=\sum_{p=0}^n q_{p,n} T^p$ on the interval $[-1, 0[$ are given by 
$$\forall \nu \in ]0,1], \;  q_n^\infty(-\nu) = \frac{1}{2\nu}\lim_{d\to +\infty}\frac{\E_{\nu}(\chi)}{(n+1)!^df_n(K)}.$$
\end{cor} 

\begin{pf}  From  Theorem~\ref{Thm_E(tilde{q})}, we get 
$$\begin{array}{rcl}
\E_{\nu}(T\tilde{q}_{V_\epsilon}(T))&=&q_K(T)-\nu q_K(\nu T)-(1-\nu)q_K((1-\nu)T)\\
&=&\frac{1}{T}[R_K(T)-R_K(\nu T)-R_K((1-\nu)T)].\\
\end{array}$$

The first part of the result follows by setting $T=-1$ and applying the property $R_K(-1+\nu)=(-1)^{n+1}R_K(-\nu)$ given by \cite{Mac} (see also \cite{SW1}).

The second part is then obtained after performing $d$ barycentric subdivisions to $K$, dividing the both sides by ${f_n(K)(n+1)!^d}$ and letting $d$ go to  $+\infty$.  Indeed, the Euler characteristic of $K$ is invariant under barycentric subdivisions, so that by \cite{DPS}, the right hand side $\frac{-2R_{\Sd^d(K)}(-\nu)}{f_n(K)(n+1)!^d}$ converges to $2\nu q_n^\infty(-\nu)$ as $d$ goes to $+\infty.$ Hence the result.
\end{pf}

\begin{cor}  \label{Cor_limsup1}Let $k=1$. For every finite  $n$-dimensional simplicial complex $K$, every $\nu\in [0,1]$ and every $0\leq i\leq n-1$,

 $$\lim_{d\to +\infty} \frac{\E(\chi)}{(n+1)!^df_n(K)}=q_n^\infty(-\frac{1}{2})= \begin{cases} 
      0 & \textup{if } n\equiv0[2] \\
      >0 &\textup{if } n\equiv 1[4] \\
      <0 &\textup{if } n\equiv 3[4]. 
   \end{cases}$$
   
   Moreover, 
$$\limsup_{d\to +\infty} \frac{\E_\nu(b_i)}{(n+1)!^df_n(K)}\leq q_{i+1,n}(1-\nu^{i+2}-(1-\nu)^{i+2})$$ and 
$$\limsup_{d\to +\infty} \sum_{j=0}^{i} (-1)^{i-j} \frac{ \E_{\nu}(b_j)}{(n+1)!^df_n(K)}\leq \sum_{j=0}^{i} (-1)^{i-j} q_{j+1,n}(1-\nu^{j+2}-(1-\nu)^{j+2}).$$
\end{cor}

Recall that  $q_{i,n}, 0\leq i\leq n$, are the coefficients of the polynomial $q_n^\infty(T)$.

\begin{pf}  
The computation $\lim_{d\to +\infty} \frac{\E(\chi)}{(n+1)!^df_n(K)}=q_n^\infty(-\frac{1}{2})$ is obtained from the second part of Corollary~\ref{Cor_E(chi)1} after performing $d$ barycentric subdivisions to $K$, dividing both sides by ${(n+1)!^df_n(K)}$ and taking the limit as $d$ tends to  $+\infty$, since the Euler characteristic of $K$ is invariant under barycentric subdivisions and $\frac{q_K}{(n+1)!^df_n(K)}$ converges to $q_n^\infty$ by \cite{DPS}. 

From \cite{DPS}, we know that the roots of $Tq_n^\infty(T)$ are symmetric with respect to $T\to -1-T$. This implies  that $q_n^\infty(-\frac{1}{2})=0$  when $n$ is even since $Tq_n^\infty(T)$ has then an odd number of roots. 
By  \cite{BW}, we know that all the roots of $Tq_n^\infty(T)$ are simple and lie in the interval $[-1, 0]$. Hence, when  $n=4s+1$  (respectively $n=4s+3$)  for some $s\in\N$, 
the polynomial $Tq_n^\infty(T)$ has $2s+1$ (respectively $2s+2$) roots between  $-\frac{1}{2}$ and 0.  Since $Tq_n^\infty(T)$ is positive for $T>0$, $q_n^\infty(-\frac{1}{2})$ is positive if there is an odd number of roots in $]-\frac{1}{2}, 0]$, negative otherwise. Hence the first part  of Corollary~\ref{Cor_limsup1}.

Now, by the second part of Theorem~\ref{Thm_E(tilde{q})}, for every $0\leq i \leq n-1,$ $$\E_\nu(\tilde{f}_i(V_\epsilon))=f_{i+1}(K)(1-\nu^{i+2}-(1-\nu)^{i+2})$$  so that performing $d$ barycentric subdivisions on $K$ we get
$$\frac{\E_\nu(b_i(V_\epsilon))}{(n+1)!^df_n(K)}\leq \frac{\E_\nu(\tilde{f}_i(V_\epsilon))}{(n+1)!^df_n(K)}=\frac{f^d_{i+1}(K)(1-\nu^{i+2}-(1-\nu)^{i+2})}{(n+1)!^df_n(K)}.$$ By \cite{DPS}, (see also\cite{SW1})  $q_{i+1,n}=\lim_{d\to +\infty}\frac{f^d_{i+1}(K)}{(n+1)!^df_n(K)}$ so that the second part follows by taking the limit. 

Finally, the last part follows along the same lines from Corollary~\ref{Cor_$E(bi)$}. Hence the result. \end{pf}

\begin{rem} 
Unlike $q_n^\infty(T)$, the polynomial $q_K(T)$ might have in general complex or non simple roots. However,  it has real coefficients and  when $K$ is a compact triangulated homology $n$-manifold, it is symmetric with respect to $T\to -1-T$, so that provided it has no  complex root with real part $-\frac{1}{2}$, we deduce as in the proof of Corollary~\ref{Cor_limsup1} the weaker equalities:

$$q_K(-\frac{1}{2})=\begin{cases} 
      0 &  \textup{if } n\equiv0[2] \\
      \geq 0 &  \textup{if } n\equiv 1[4] \\
      \leq 0 &  \textup{if } n\equiv 3[4].
   \end{cases}$$ 
\end{rem}

We finally deduce a probabilistic proof of the following theorem.

\begin{thm} \label{Thm_Euler$k=1$} If $K$ is a compact triangulated homology $2n$-manifold, then 
$$\chi(K)=\sum\limits_{p=0}^{2n}(-\frac{1}{2})^p  f_p(K). $$
\end{thm}

Theorem~\ref{Thm_Euler$k=1$} has already been proved in \cite{A}, see also \cite{SW1}. 
 
\begin{pf} It follows from Corollary~\ref{Cor_VHom} that for every $\epsilon\in \mC^0(K)$, the hypersurface $V_\epsilon$ is the triangulation of an odd-dimensional homology
manifold, so that $\chi(V_\epsilon)$ vanishes from Poincar\'e duality, see \cite{M}. We deduce that $\E(\chi)$ vanishes as well and so the result follows from Corollary~\ref{Cor_E(chi)1}.   \end{pf}

\begin{rem} \label{Rem_zerosqinfity}
Note that Theorem~\ref{Thm_Euler$k=1$}  implies that  $q_n^\infty(-\frac{1}{2})=0$ when $n>0$ is even and $q_n^\infty(-{1})=0$ for every $n>0$. This fact can be deduced from the symmetry property of  $R_K(T)=Tq_K(T)-\chi(K)T$ see \cite{Mac} and also \cite{SW1}.

Note also that from Corollary~\ref{Cor_E(chi)1} we more generally deduce that under the hypothesis of 
Theorem~\ref{Thm_Euler$k=1$}, $\chi(K)=\nu q_K(-\nu)+(1-\nu)q_K(\nu-1)$ for every $\nu\in [0,1].$ This formula has been observed by 
A.~Kassel right after the first author  gave an informal talk for students on Theorem~\ref{Thm_Euler$k=1$}.

\end{rem}

\subsection{The case of compact homology manifolds}\label{Ssect_HomMfd}
When $K$ is a triangulated compact homology manifold, it inherits a decomposition into blocks which is dual  to the triangulation and is useful to prove Poincar\'e duality, see \cite{M}. In particular, these blocks span a chain complex
which computes the homology of $K$ exactly as if it were a $\CW$-complex (from the homology point of view, there is no difference).

By definition, the block $D(\sigma)$ dual to a simplex $\sigma\in K$ is the union of all open simplices $[\hat{\sigma}_0\ldots, \hat{\sigma}_p]$ of $\Sd(K)$ such that $\sigma_0=\sigma$. The union of closed such simplices is denoted by $\overline{D}(\sigma).$

\begin{lem} \label{Lem_Vdual}
Let $K$ be a triangulated homology manifold of dimension $n>0$. Then, for every $1\leq k\leq n$ and every $\epsilon\in \mC^{k-1}(K)$, $V_\epsilon$ is the union of the blocks $D(\sigma)$  dual to the  simplices $\sigma$ of $K$ such that the restriction  of $d\epsilon$ to $\sigma$ does not vanish.
\end{lem}

\begin{pf}
This follows  from the definitions of $V_\epsilon$ and the dual block decompositions of $K$.
\end{pf}

For every $\epsilon \in \mC^{k-1}(K)$, every $i\in \{0,\ldots, n-k\}$ and every $\nu\in [0,1]$, let us denote by $\check{f}_i(V_\epsilon)$  the number of $i$-dimensional blocks of $K$ that are in $V_\epsilon$ and by $\E_\nu(\check{f}_i)=\int_{\epsilon\in \mC^{k-1}(K)}\check{f}_i(V_\epsilon)d\mu_{\nu}(\epsilon)$  its  mathematical expectation.

\begin{thm}\label{Thm_vfi} 
Let $K$ be a compact triangulated homology manifold of dimension $n>0$. Then, for every $1\leq k\leq n$,  every $0\leq i\leq n-k$ and every $\nu\in [0,1],$ $\E_\nu(\check{f}_i)=m_k(K^{[n-i]})$, so that $\E_\nu(\check{q}_{V_\epsilon}(T))=\int_K\check{q}_{D(\sigma)}(T)dm_k(\sigma)$.
\end{thm}

Note that for every $\sigma\in K,$ $D(\sigma)$ is made of a single block, so that $\check{q}_{D(\sigma)}(T)=T^{n-\dim\sigma}.$
\begin{pf}
By definition,
$$\begin{array}{lcl}
\E_\nu(\check{f}_i)&=&\int\limits_{\epsilon \in \mC^{k-1}(K)}  \sum\limits_{\tiny{\begin{array}{c}\sigma\in K^{[n-i]}\\ D(\sigma)\in V_\epsilon\end{array}}}1 d\mu_\nu(\epsilon)\\
&&\\
&=&  \sum\limits_{\tiny{\sigma\in K^{[n-i]}}} \mu_\nu\{\epsilon \in \mC^{k-1}(K) |\hat{\sigma}\in V_\epsilon\}.\\
\end{array}$$
We thus deduce from Lemma~\ref{Lem_Vdual} and the definition of $m_k$ that
$$\E_\nu(\check{f_i})=m_k(K^{[n-i]}).$$
Summing over all $i\in\{0,\ldots,n-k\}$ we get $\E_v(\check{q}_{V_\epsilon}(T))=\int_K\check{q}_{D(\sigma)}dm_k(\sigma),$ since by definition $\check{q}_{D(\sigma)}(T)=T^{n-\dim\sigma}$.
\end{pf}

\begin{cor}\label{Cor_vbi}
Under the hypotheses of Theorem~\ref{Thm_vfi}, the following average Morse inequalities hold true:
$$\E_\nu(b_i)\leq   m_k(K^{[n-i]})$$ 
and 

$$\sum_{j=0}^i(-1)^{i-j}\E_\nu(b_j)\leq \sum_{j=0}^{i}(-1)^{i-j}m_k(K^{[n-j]}),$$ where the latter is an equality if $i=n-k$.
\end{cor}

\begin{pf} It follows from Lemma~\ref{Lem_Vdual} that for every $\epsilon\in \mC^{k-1}(K)$ the blocks in $V_\epsilon$ provide a filtration of $V_\epsilon$ so that the hypotheses of Theorem~39.5 and Theorem~64.1 of \cite{M} are satisfied.
As a consequence, the chain complex spanned by the blocks of $V_\epsilon$ compute the homology of $V_\epsilon$ with $\Z/2\Z$-coefficients and thus the result follows from Theorem~\ref{Thm_vfi}  and the Morse inequalities.
\end{pf}

\begin{cor} Let $K$ be a compact triangulated homology manifold of dimension $n>0$. Then, for every $1\leq k\leq n$ and  every $\nu\in [0,1],$
$$\E_\nu(\chi)=\int_K(-1)^{n-\dim\sigma}dm_k(\sigma).$$

In particular,

$$(-1)^{n+1}\E_\nu(\chi)+\chi(K)=\sum\limits_{p=0}^{n}(-1)^p\mu_\nu(Z^{k-1}(\Delta_p))f_p(K),$$
where $\mu_\nu(Z^{k-1}(\Delta_p))=1$, if $p<k$.
\end{cor}

\begin{pf}
As in Corollary~\ref{Cor_$E(bi)$}, the case of equality $i=n-k$ in Corollary~\ref{Cor_vbi} computes the expected Euler characteristic of $V_\epsilon$, hence the first part. By Proposition~\ref{Prop_mk}, we then deduce $(-1)^n\E_\nu(\chi)=\sum_{p=k}^{n}(-1)^p\big(1-\mu_\nu(Z^{k-1}(\Delta_p))\big)f_p(K)$ and the second part. 
\end{pf}

\begin{rem} Corollary~\ref{Cor_vbi} holds true for every $k$ provided $K$ is a compact triangulated homology manifold while 
Corollary~\ref{Cor_$E(bi)$} holds true for every finite simplicial complex provided $k=1$. When both conditions are satisfied, that is $k=1$ and $K$ is a triangulated compact homology $n$-manifold, then from Corollary~\ref{Cor_VHom}, $V_\epsilon$ is a compact triangulated homology $(n-1)$-manifold for every $\epsilon\in \mC^{k-1}(K)$ so that  Poincar\'e duality holds true for $V_\epsilon$ by Theorem~65.1 of \cite{M}. Then, for every $0\leq i\leq n-1$, $\E_{\nu}(b_i)=\E_{\nu}(b_{n-1-i})$, so that the upper estimates of Corollary~\ref{Cor_$E(bi)$}  and Corollary~\ref{Cor_vbi}   coincide after the change  $i\mapsto n-1-i.$ 

\end{rem}

\begin{cor} Under the hypotheses of Theorem~\ref{Thm_vfi},

$$\limsup_{d\to +\infty}\frac{\E_\nu(b_i)}{(n+1)!^df_n(K)}\leq q_{n-i,n}\big(1-\mu_{\nu}(Z^{k-1}(\Delta_{n-i}))\big)$$ and 

$$\limsup_{d\to +\infty} \sum_{j=0}^i\frac{(-1)^{i-j}\E_\nu(b_i)}{(n+1)!^df_n(K)}\leq  \sum_{j=0}^i (-1)^{i-j}q_{n-j,n}\big(1-\mu_{\nu}(Z^{k-1}(\Delta_{n-j}))\big)$$  where the latter is an equality if $i=n-k$.
\end{cor}

\begin{pf}
It follows from Corollary~\ref{Cor_vbi}, along the same lines as Corollary~\ref{Cor_limsup1}. 
\end{pf}

\subsection{ The general case} \label{Ssect_$k>1$} 
When $K$ is a general finite simplicial complex, the dual block decomposition of $K$ does not span a chain complex which computes the homology of $K$ so that the results of Section~\ref{Ssect_HomMfd}  do not apply. Likewise, Proposition~\ref{Prop_$k=1$}
does not extend to $k>1$, so that in this case the results of Section~\ref{Ssect_$k=1$} do not apply. However, for every $\epsilon \in \mC^{k-1}(K),$ $V_\epsilon$ is a subcomplex  of $\Sd(K)$ so that its homology can be computed with the help of the simplicial homology theory, providing  weaker upper estimates which we are going to obtain now.

Recall that every simplex $\alpha\in \Sd(K)$ is of the form $\alpha=[\hat{\sigma}_0,\ldots, \hat{\sigma}_p],$ where $\sigma_0<\ldots<\sigma_p$ are simplices of $K$. We set, following \cite{M}, $\ini(\alpha)=\sigma_0$ and $\fin(\alpha)=\sigma_p$.

\begin{thm} \label{Thm_Enu} Let $K$ be a finite $n$-dimensional simplicial complex and $k\in\{ 1,\ldots, n\}$. Then, for every $\nu\in[0,1]$, $$\E_{\nu}(q_{V_\epsilon}(T))=\int_Kq_{D(\sigma)}(T)dm_k(\sigma).$$

Moreover, If $K$ is a compact triangulated homology $n$-manifold and $k=1$, then $$\E_\nu(R_{V_\epsilon}(-1-T))=(-1)^{n} \E_\nu({R_{V_\epsilon}(T)}),$$ where $R_{V_\epsilon}(T)=Tq_{V_\epsilon}(T)-\chi(V_\epsilon)T.$
\end{thm}

\begin{pf} We first observe that for every $0\leq i\leq n-k$ and every $\nu\in [0,1],$
$$
\begin{array}{lcl}
\E_\nu(f_i)&=&\sum\limits_{\epsilon\in C^{k-1}(K)}  \sum\limits_{\tiny{\begin{array}{c}\alpha\in \Sd(K)^{[i]}\\ \alpha \in V_\epsilon\end{array}}}\mu_\nu(\epsilon)\\
&& \\
&=& \sum\limits_{\alpha\in \Sd(K)^{[i]}} \mu_\nu\{\epsilon\in C^{k-1}(K)| \alpha\in V_\epsilon\}.
\end{array}
$$

We then deduce
$$
\begin{array}{lcl}
\E_\nu(q_{V_\epsilon}(T))&=&\sum\limits_{i=0}^{n-k} T^i \sum\limits_{\alpha\in \Sd(K)^{[i]}}\mu_\nu\{\epsilon\in C^{k-1}(K)| \alpha\in V_\epsilon\}\\
&& \\
&=& \sum\limits_{\sigma\in K}\sum\limits_{\tiny{\begin{array}{c}\alpha\in \Sd(K)\\ \ini(\alpha)=\sigma\end{array}}} T^{\dim\alpha}\
 \mu_\nu\{\epsilon\in C^{k-1}(K)| \alpha\in V_\epsilon\}\\
 && \\
 &=&\int_Kq_{D(\sigma)}(T)dm_k(\sigma).
\end{array}
$$

Finally, when $K$ is a compact triangulated homology manifold and $k=1$, we know from Corollary~\ref{Cor_VHom} that for every $\epsilon\in \mC^{k-1}(K)$,  $V_\epsilon$ is itself a compact triangulated homology manifold. 
Then, from Theorem~2.1 of \cite{Mac} (see also \cite{SW1}), it follows that $R_{V_\epsilon}(-1-T)=(-1)^{n}R_{V_\epsilon}(T).$
The result thus follows after integration over $\mC^{k-1}(K).$
\end{pf}

\begin{cor} \label{Cor_Enuchi}Let $K$ be a finite $n$-dimensional simplicial complex and $k\in\{1,\ldots, n\}$, $\nu\in[0,1]$. Then, $$\E_{\nu}(\chi)=\int_K\left(\chi(D(\sigma))-\chi(\dot{D}(\sigma))\right) dm_k(\sigma).$$

Moreover, for every $0\leq i\leq n-k$, the following average Morse inequalities hold:
$$\E_{\nu}(b_i)\leq \E_{\nu}(f_i)=\int_Kf_i(D(\sigma))dm_k(\sigma)$$ and 
$$\sum_{j=0}^i(-1)^{i-j}\E_{\nu}(b_j)\leq \int_K\sum_{j=0}^{i}(-1)^{i-j}f_j(D(\sigma))dm_k(\sigma),$$
where the latter is an equality if $i=n-k$.
\end{cor}

\begin{pf} The first part follows from Theorem~\ref{Thm_Enu}  after evaluation at $T=-1$, since by definition, for every $\epsilon\in \mC^{k-1}(K),$ $\chi(V_\epsilon)=q_{V_\epsilon}(-1)$. Then, for every $\epsilon\in \mC^{k-1}(K)$, the Morse inequalities applied to the simplicial  chain complex of $V_\epsilon$ with $\Z/2\Z$-coefficients read $b_i(V_\epsilon)\leq f_i(V_\epsilon)$ and $\sum_{j=0}^i (-1)^{i-j}b_j(V_\epsilon)\leq \sum_{j=0}^i(-1)^{i-j}f_j(V_\epsilon)$, the latter being an equality when $i=n-k$. The last part of Corollary~\ref{Cor_Enuchi} thus follows after integration over $ \mC^{k-1}(K).$
\end{pf}

Let us denote by $\lambda_{i,j}$ the number of interior $(j-1)$-faces of the subdivided standard simplex $\Sd(\Delta_{i-1}),$ with the convention that $\lambda_{i,0}=0$ if $i>0$ and $\lambda_{0,0}=1$.

\begin{defn} \label{Defn_ $c_i^+(n,k)$}
 For every $0\leq i\leq  n-k$, we set $$c_i^+(n,k)=\sum_{p=k+i}^n\delta^{p,k}_i q_{p,n}$$
  where 
  $\delta_i^{p,k}=\sum\limits_{l=k}^{p-i}\big(1-\mu_\nu(Z^{k-1}(\Delta_l))\big)f_l(\Delta_p)\lambda_{p-l,i}.$
 \end{defn}

The constant $\delta_{i}^{p,k}$ appearing in Definition~\ref{Defn_ $c_i^+(n,k)$} has actually a probabilistic interpretation, namely $\delta_i^{p,k}=\int_{C^{k-1}(\Delta_p)}f_i(V_\epsilon\cap \stackrel{\circ}{\Delta}_p)d\mu_\nu(\epsilon).$
We develop this interpretation in Appendix A (see (\ref{Eqn_delta})).

\begin{cor} \label{Cor_EbiEuler$k>1$} Let $K$ be a finite $n$-dimensional simplicial complex. Then, for every $k\in\{1,\ldots, n\}$ and every $0\leq i \leq n-k$,

$$\limsup_{d\to +\infty}   \frac{\E_\nu(b_i)}{(n+1)!^df_n(K)}\leq c_i^+(n,k), $$ 
$$\limsup_{d\to +\infty} \frac{1}{(n+1)!^df_n(K)}(\sum\limits_{j=0}^{i}(-1)^{i-j}\E_\nu(b_j))\leq \sum\limits_{j=0}^{i}(-1)^{i-j}c_j^+(n,k)$$
and

$$\lim_{d\to +\infty}   \frac{\E_\nu(\chi)}{(n+1)!^df_n(K)}=  \sum\limits_{i=0}^{n-k}(-1)^ic^+_i(n,k).$$
\end{cor}

\begin{pf} By Proposition~\ref{Prop_mk}, $m_k=\sum_{p=k}^{n}\big(1-\mu_\nu(Z^{k-1}(\Delta_p))\big)\gamma_{p,K}$ while  $\gamma_{p,K}^d=\frac{\gamma_{p, \Sd^d(K)}}{(n+1)!^d}$ weakly converges to $dvol_K$ as $d$ grows to $+\infty$, see \cite{SW1}.
By Corollary~\ref{Cor_Enuchi}, for every $0\leq i\leq n-k$, $\frac{\E_{\nu}(b_i)}{(n+1)!^df_n(K)}=\frac{1}{f_n(K)}\sum_{p=k}^{n}\big(1-\mu_\nu(Z^{k-1}(\Delta_p))\big)\int_K f_i(D(\sigma))d\gamma_{p,K}^d(\sigma)$.
Now, by  \cite{SW1}, $\frac{1}{f_n(K)}\int_Kf_i(D(\sigma))d\gamma_{p,K}^d(\sigma)$ converges to $\sum_{h=i}^{n-p}q_{p+h,n}f_p(\Delta_{p+h})\lambda_{h,i}$ as $d$ grows to $+\infty$.
We deduce that  $$\limsup_{d\to +\infty}\frac{\E_\nu(b_i)}{(n+1)!^df_n(K)}\leq c_i^+(n,k)$$ with 
$$\begin{array}{lcl}
c_i^+(n,k)&=&\sum_{p=k}^{n-i}\big( 1-\mu_\nu(Z^{k-1}(\Delta_p))\big)\sum_{h=p+i}^nq_{h,n}f_p(\Delta_h)\lambda_{h-p, i}\\
&=&\sum_{h=k+i}^{n}q_{h,n}\sum_{p=k}^{h-i}\big(1-\mu_\nu(Z^{k-1}(\Delta_p))\big)f_p(\Delta_h)\lambda_{h-p,i}.\end{array}$$
Hence the first part of the result.

The second part just follows the from the Morse inequalities in Corollary~\ref{Cor_Enuchi}.
As for the last part, it follows from Corollary~\ref{Cor_Enuchi} and  what we have just done, since $\E_\nu(\chi)=\sum_{i=0}^{n-k}(-1)^i\int_K f_i(D(\sigma))dm_k(\sigma).$
\end{pf}

\section{Asymptotic topology of $V_\epsilon$, lower estimates} \label{Sect_Lower}
\subsection{Lower estimates for the expected Betti numbers} 
Let us  start with a key proposition, the proof of which is inspired by  H.~Whitney's proof of the existence of triangulation on smooth manifolds~\cite{W}.

\begin{prop} \label{Prop_Pair} 
For every closed (not necessarily connected) submanifold   $\Sigma\subset \R^n$  of codimension $k\geq1$,
there exists $m_n>0$ such that for every $ m\geq m_n$  the pair $(\stackrel{\circ}{\Delta}_n, V_\epsilon)$ gets homeomorphic to $(\R^n, \Sigma)$ for some $\epsilon \in \mC^{k-1}(\Sd^m(\Delta_n))$, where $\Delta_n$ denotes the standard $n$-simplex. 

\end{prop}
In the light of Proposition~\ref{Prop_Pair} , let us set the following.
\begin{defn}\label{Defn_Complexity} 
 The {complexity} of a closed (not necessarily connected) submanifold $\Sigma$ of $\R^n$ is the smallest value $m_n$ given by Proposition~\ref{Prop_Pair}. It is denoted by $m_n(\Sigma)$. Likewise, the {$n$-dimensional complexity} of a closed connected manifold $\Sigma$  which embeds into  $\R^n$ is the infimum of $m_n(\Sigma)$ over all embeddings $\Sigma\hookrightarrow\R^n$.
\end{defn}

\begin{rem}\label{Rem_Complexity} 
  Recall that from H.~Whitney's embedding theorem, every  manifold  $\Sigma$ of dimension $n$ embeds in $\R^{2n+1}$. 
Proposition~\ref{Prop_Pair} and Definition~\ref{Defn_Complexity} thus provide a combinatorial complexity of such closed $n$-dimensional manifold, namely  the infimum of $m_{2n+1}( \Sigma)$ over all embeddings $\Sigma\hookrightarrow \R^{2n+1}$.
\end{rem}

\textbf{Proof of Proposition~\ref{Prop_Pair}.}  
Consider a diffeomorphism $\phi: \R^n\to \stackrel{\circ}{\Delta}_n$ and set $\widetilde{\Sigma}=\phi(\Sigma)$.
There is a  positive integer $m_n$ such that for every $m\geq m_n$, possibly after a small perturbation of $\phi$ by an isotopy of $\stackrel{\circ}{\Delta}_n$ with compact support, we  have  the following properties for $\widetilde{\Sigma}$ (conditions that appear in Section 13 of~\cite{W}) :

\begin{enumerate}
\item $\widetilde{\Sigma}$ does not intersect the $(k-1)$-skeleton of $\Sd^m(\Delta_n)$,
\item the intersection, ${\widehat{\Sigma}_k}$, of $\widetilde{\Sigma}$ with each $k$-simplex of $\Sd^m(\Delta_n)$, if  not empty, is transversal at one point,
\item the intersection of $\widetilde{\Sigma}$ with each $l$-simplex  $\sigma$ of dimension $l>k$ is isotopic to the cone over the intersection of ${\widehat{\Sigma}_{l-1}}$ with the $(l-1)$-skeleton of $\sigma$ centered at the barycenter $\hat{\sigma}$ of $\sigma$.
\end{enumerate}

Let $\widehat{\Sigma}=\widehat{\Sigma}_n$ denote the cone described above. The pair $(\stackrel{\circ}{\Delta}_n, \widetilde{\Sigma})$ is homeomorphic to $(\stackrel{\circ}{\Delta}_n, \widehat{\Sigma})$.

Let $\hat{\epsilon}\in \mC^k(\Sd^m(\Delta_n))$ be the $k$-cochain defined by the relation  $<\hat{\epsilon}, \sigma>=\widehat{\Sigma}\circ \sigma$ which is  either $0$ or  $1$ for every  $k$-simplex $\sigma \in \Sd^m(\Delta_n)$ (here $\circ$ denotes the intersection index of $\Sigma$ and $\sigma$). By definition, $\hat{\epsilon}$ is a cocycle. Indeed, for all $(k+1)$-simplex $\tau$ of $\Sd^m(\Delta_n)$, the intersection of $\widehat{\Sigma}$ and  $\tau$ is either empty or isotopic to an interval. Thus, we have $<d \hat{\epsilon}, \tau>=<\hat{\epsilon}, \partial \tau>=0$ for all $(k+1)$-simplex $\tau\in  \Sd^m(\Delta_n)$. 
Therefore, as the simplex $\Delta_n$ is acyclic, there exists $\epsilon \in \mC^{k-1}(\Sd^m(\Delta_n))$ such that $\hat{\epsilon}=d\epsilon$. 

To construct the isotopy between $V_\epsilon$ and $\widehat{\Sigma}$, we proceed by induction on the dimension of the skeleton as in the proof of Proposition~\ref{Prop_$k=1$}. Let us recall that the intersection of a $k$-simplex $\sigma$  with $V_\epsilon$ is either empty or a single point, as in the case of the intersection of $\sigma$ 
with  $\widehat{\Sigma}$. Moreover,  by definition of $\epsilon$, $V_\epsilon \cap \sigma$ is non empty if and only if $\widehat{\Sigma}\cap \sigma$ is non empty. When the intersection is non-empty, we consider the two points on $\sigma$ determined as the intersection with $V_\epsilon$ and with $\widehat{\Sigma}$ and isotope $V_\epsilon$ so that these two points match. Now let us suppose that $V_\epsilon$ can be isotoped to $\widehat{\Sigma}$ up to the level of $l$-skeleton for $l\geq k$. Let $\tau$ be a $(l+1)$-simplex.
If $V_\epsilon$ meets $\tau$, then the intersection $V_\epsilon \cap \tau$ is the cone centered at the barycenter $\hat{\tau}$ of $\tau$ over the intersections of $V_\epsilon$ with the $l$-skeleton of $\tau$. Besides $\widehat{\Sigma}\cap \tau$, if not empty, is the cone defined at $\hat{\tau}$ over the intersection of $\hat{\Sigma}$ with the $l$-skeleton of $\tau$. 

We can extend the isotopy between  $\Sigma\cap \partial \tau$ and $V_\epsilon \cap \partial \tau$ to an isotopy between $\Sigma\cap \tau$ and $V_\epsilon \cap \tau$ by taking the cone over $\hat{\tau}$.
\hfill $\square$

\begin{rem} \label{Rem_pair} 
 
Note that  the cochain $\epsilon\in \mC^{k-1}(\Sd^m(\Delta_n))$ given by Proposition~\ref{Prop_Pair}  is not unique since we may add to $\epsilon$ any $(k-1)$-cocycle.
For example, when $k$ is odd, another $\epsilon$ is obtained by switching the labelings of 0 and 1, that is replacing $\epsilon$ by $1-\epsilon$. When $k$ is even, consider for example the cochain  $\tilde{\epsilon}\in \mC^{k-2}(\Sd^m(\Delta_n))$  which labels 1 only one $(k-2)$-dimensional simplex. Then $d\tilde{\epsilon}$ is a $(k-1)$-cocycle and  we may replace $\epsilon$ by $\epsilon+d\tilde{\epsilon}.$ 
\end{rem}

Let $\Sigma\subset \R^n$ and $m_n(\Sigma)$ be given by Proposition~\ref{Prop_Pair}. Let $m_n(\Sigma)\leq m\leq d$ and  $\sigma\in \Sd^{d-m}(K)$ be an $n$-simplex. We set $$prob(\sigma, \Sigma)=\mu_\nu\{\epsilon \in \mC^{k-1}(\Sd^d (K))| (\stackrel{\circ}{\sigma}, V_\epsilon\cap\stackrel{\circ}{\sigma}) \mbox{ is homeomorphic to } (\R^n, \Sigma)\}.$$

\begin{thm}\label{Thm_$p_Sigma$} For every $n$-dimensional finite simplicial complex $K$, every closed codimension $k\geq 1$ submanifold $\Sigma$ of $\R^n$, every  $m_n(\Sigma)\leq m\leq d$  and every $n$-simplex $\sigma\in \Sd^{d-m}(K)$,   $$prob(\sigma, \Sigma) \geq \frac{1}{2^{f_{k-1}^m(\Delta_n)-1}}$$ where  $f_{k-1}^m(\Delta_n)=f_{k-1}(\Sd^m(\Delta_n)).$ 
\end{thm}

\begin{pf}
Observe that $f_{k-1}^{d}(K)=f_{k-1}^{m}(\sigma)+ R$  where  $R$ is the number of $(k-1)$-simplices of $\Sd^d(K)$ which are not in $\Sd^m(\sigma)$. 
By Remark~\ref{Rem_pair},  we know that there are at least two choices of ${\epsilon_\sigma}\in \mC^{k-1}(\Sd^m(\sigma))$ with the property that the pair $(\stackrel{\circ}{\sigma}, V_{{\epsilon_\sigma}})$ is homeomorphic to $(\R^n, \Sigma)$.
Any of these two $\epsilon_\sigma$ extend to $\epsilon \in \mC^{k-1}(\Sd^d(K))$ such that $(\stackrel{\circ}{\sigma}, V_\epsilon\cap \stackrel{\circ}{\sigma})$ is homeomorphic to $(\R^n, \Sigma)$.
The restriction of $\epsilon$ on the $R$ many $(k-1)$-simplices of $\Sd^d(K)\setminus \Sd^m(\sigma)$ is arbitrary and so we deduce that 
$prob(\sigma,\Sigma)$  is at least  $\frac{2\times2^{R}}{2^{f_{k-1}^d(K)}}=\frac{2\times2^{R}}{2^{f_{k-1}^{m}(\sigma)+R}}=\frac{1}{2^{f_{k-1}^{m}(\Delta_n)-1}}$.
\end{pf}

For every $\epsilon\in \mC^{k-1}(\Sd^{d}(K))$, let  $N_{\Sigma} (\epsilon)$ be the maximum number of disjoint open simplices $(\stackrel{\circ}\sigma_j)_{j\in J}$ which can be packed in $K$ in such a way that  the pair $(\stackrel{\circ}\sigma_j, V_\epsilon \cap \stackrel{\circ}\sigma_j)$ is homeomorphic to $(\R^n, \Sigma)$, where for every $j\in J$,  $\sigma_j$ is a simplex of $\Sd^{d-m_j}(K)$ for some $m_j\in \{{m_n(\Sigma)},\ldots, d\}$. 
We now set  $p_\Sigma=\frac{1}{2^{f_{k-1}^{m_n(\Sigma)}(\Delta_n)-1}},$ the right hand side of the inequality in Theorem~\ref{Thm_$p_Sigma$} for $m=m_n(\Sigma)$ and  $c_\Sigma= \frac{p_\Sigma}{(n+1)!^{m_n(\Sigma)}}>0$ and $\E_\nu(N_{\Sigma})=\int_{\mC^{k-1}(\Sd^d(K))} N_{\Sigma}(\epsilon)d\mu_\nu(\epsilon)$.

\begin{thm}  \label{Thm_cSigma} For every  $n$-dimensional finite simplicial complex $K$   and every closed codimension $k\geq 1$ submanifold $\Sigma$ of $\R^n$,  \[c_\Sigma\leq \liminf_{d\to +\infty} \frac{\E_\nu(N_{\Sigma})}{ (n+1)!^df_n(K)}. \]
\end{thm}

\begin{pf} Let $d\geq m_n(\Sigma)$. For every n-simplex $\sigma\in \Sd^{d-m_n(\Sigma)}(K)$ and every $\epsilon\in \mC^{k-1}(\Sd^{d}(K))$,
let $N_{\Sigma, \sigma}(\epsilon)$ be equal to 1 if $( \stackrel{\circ}\sigma,  V_\epsilon \cap \stackrel{\circ}\sigma)$ is homeomorphic to $(\R^n, \Sigma)$ and 0 otherwise.  Then

$$\begin{array}{lcl}
\E_\nu(N_{\Sigma})&=&\int_{ \mC^{k-1}(\Sd^d(K))} N_{\Sigma}(\epsilon) d\mu_\nu(\epsilon)\\
&\geq &\int_{ \mC^{k-1}(\Sd^d(K))} \sum\limits_{\sigma \in \Sd^{d-m_n(\Sigma)}(K)}N_{\Sigma,\sigma}(\epsilon) d\mu_\nu(\epsilon) \\
&= &\sum\limits_{\sigma \in  \Sd^{d-m_n(\Sigma)}(K)}  \int_{ \mC^{k-1}(\Sd^d(K))}N_{\Sigma,\sigma}(\epsilon) d\mu_\nu(\epsilon)\\
&\geq& (n+1)!^{d-m_n(\Sigma)} f_n(K) p_\Sigma.
\end{array}
$$

The last line follows from Theorem~\ref{Thm_$p_Sigma$} and the fact that the number of $n$-simplices  of an $n$-dimensional simplicial complex gets multiplied by $\lambda_{n+1, n+1}=(n+1)!$ after a barycentric subdivision, see \cite{DPS} or also \cite{SW1}.
Thus, we get
$$\frac{\E_\nu(N_{\Sigma})}{ (n+1)!^df_n(K)} \geq \frac{p_\Sigma}{{(n+1)!}^{m_n(\Sigma)}}=c_\Sigma.$$

Hence the result.
\end{pf}

Now, for every $m\geq 1$, let $\mmC(m)$ be the finite set of homeomorphism classes of pairs $(\R^n, \Sigma)$, where $\Sigma$ is a closed connected  $(n-k)$-dimensional manifold embedded in $\R^n$ by an embedding of complexity $m=m_n(\Sigma)$, see Definition~\ref{Defn_Complexity}. 
For every $\epsilon \in \mC^{k-1}(\Sd^d (K))$ and every $i\in \{0,1,\ldots,n-k\}$,  we set $b_i(V_\epsilon)=\dim H_i(V_\epsilon)$ and $\E_\nu(b_i)=\int_{ \mC^{k-1}(\Sd^d (K))}b_i(V_\epsilon)d\mu_\nu(\epsilon).$

\begin{defn} \label{Defn_$c_i^-$} For every $1\leq k\leq n$ and $0\leq i\leq n-k$, we set $$c_i^-(n,k)=\sum\limits_{m=1}^{+\infty} \frac{1}{2^{{f_{k-1}^m(\Delta_n)}-1}(n+1)!^m}   \sum\limits_{(\R^n, \Sigma)\in \mmC(m)} b_i(\Sigma).$$

\end{defn}

\begin{cor} \label{Cor_Inf} For every finite $n$-dimensional complex $K$, every $1\leq k\leq n$ and every $i\in\{0,1,\ldots, n-k\}$,

 \[ \liminf_{d\to +\infty} \frac{\E_\nu(b_i)}{ (n+1)!^df_n(K)}\geq c_i^-(n,k) \]
\end{cor}

  \begin{pf}
  Let $\epsilon\in \mC^{k-1}(\Sd^d (K))$. For every connected component of $V_\epsilon$ which is contained in the interior of an $n$-simplex $\sigma\in \Sd^{d-m}(K)$ in such a way that $( \stackrel{\circ}\sigma,  V_\epsilon \cap \stackrel{\circ}\sigma)$ is homeomorphic to $(\R^n, \Sigma)$ for some codimension $k$ submanifold $\Sigma$ of $\R^n$, where $m_n(\Sigma)\leq m\leq d$, the homeomorphism type of the pair $(\R^n, \Sigma)$ does not depend on the choice of $\sigma$ and $m$ in the case it is not unique. We deduce that for every $M>0$,
   $$b_i(V_\epsilon)\geq \sum\limits_{m=1}^{M} \sum\limits_{(\R^n,\Sigma)\in \mmC(m)}b_i(\Sigma)N_{\Sigma}(\epsilon).$$

After integration we get, $\E_\nu(b_i)\geq \sum\limits_{m=1}^{M}\sum\limits_{(\R^n,\Sigma)\in \mmC(m)}b_i(\Sigma)\E_\nu(N_{\Sigma}).$
Theorem~\ref{Thm_cSigma} then implies that

$$\liminf_{d\to +\infty}\frac{\E_\nu(b_i)}{(n+1)!^df_n(K)}\geq \sum\limits_{m=1}^M \sum\limits_{(\R^n,\Sigma)\in\mmC(m)} \frac{b_i(\Sigma)}{2^{f^m_{k-1}(\Delta_n)-1}(n+1)!^m}.$$
The result follows by letting $M$ grow to $+\infty$.
 \end{pf}

\subsection{Complexity of surfaces in $\R^3$}\label{Sect_Complexity}
\begin{figure}[h]
   \begin{center}
    \includegraphics[scale=0.4]{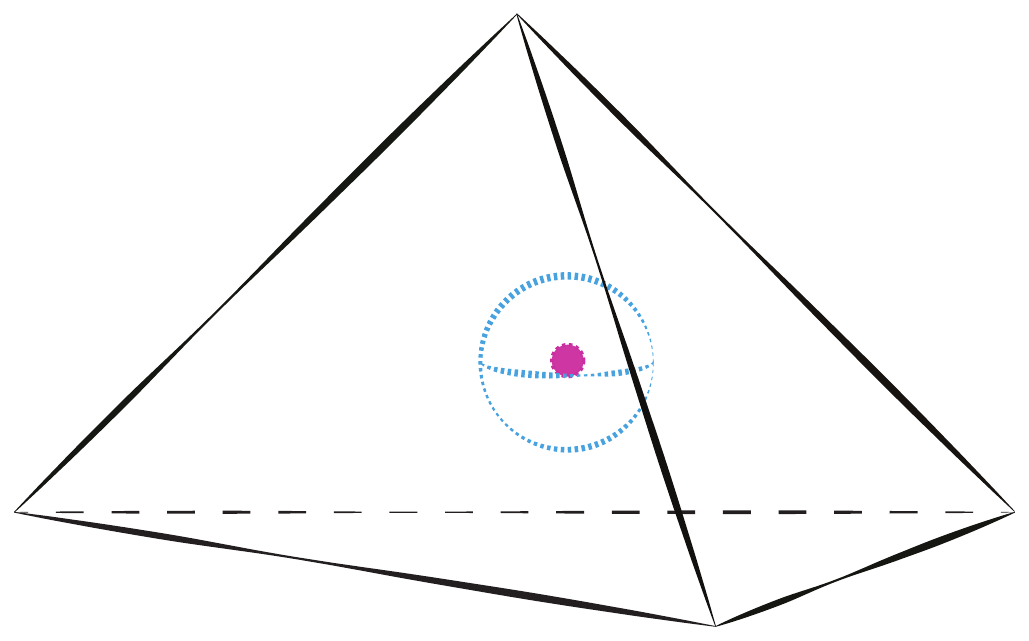}
       \caption{In $\Sd(\Delta_3)$,  $V_\epsilon$  may be a topological sphere.}
    \label{Complexity}
      \end{center}
 \end{figure}

Let us now study the 3-dimensional complexity of surfaces in the sense of Definition~\ref{Defn_Complexity}. 
We first observe that  there exists $\epsilon\in \mC^0(\Sd(\Delta_3))$  such that $V_\epsilon$ is homeomorphic to a 2-sphere. Indeed, let $\epsilon$ take the value 0 on the barycenter of $\Delta_3$ and 1 on all the other vertices of $\Sd(\Delta_3)$, see Figure~\ref{Complexity}. The complexity of the 2-sphere is thus 1.

\begin{figure}[h]
   \begin{center}
   \includegraphics[scale=0.4]{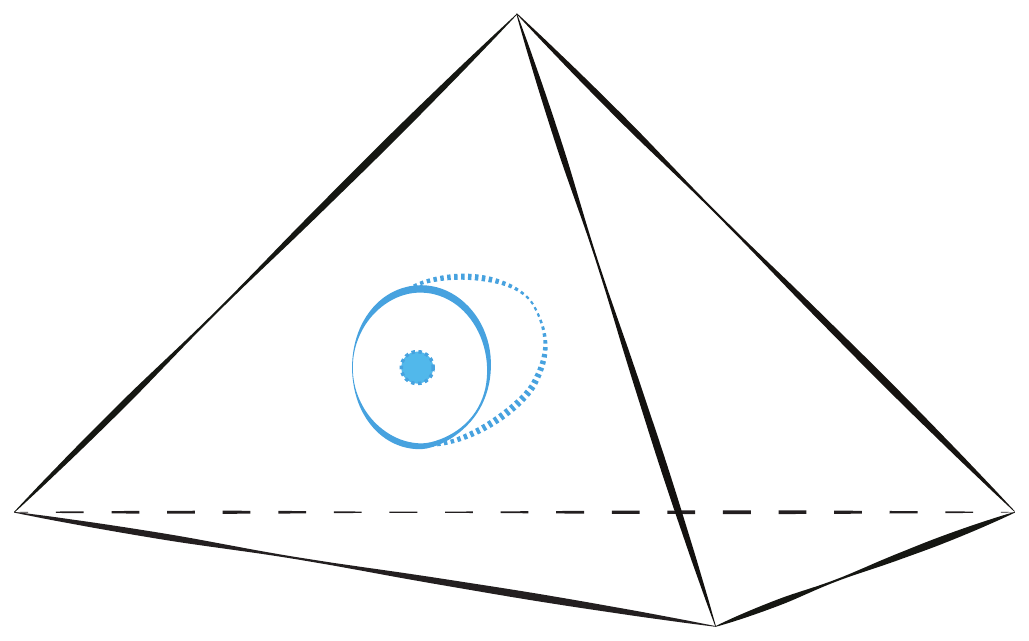}
     \includegraphics[scale=0.4]{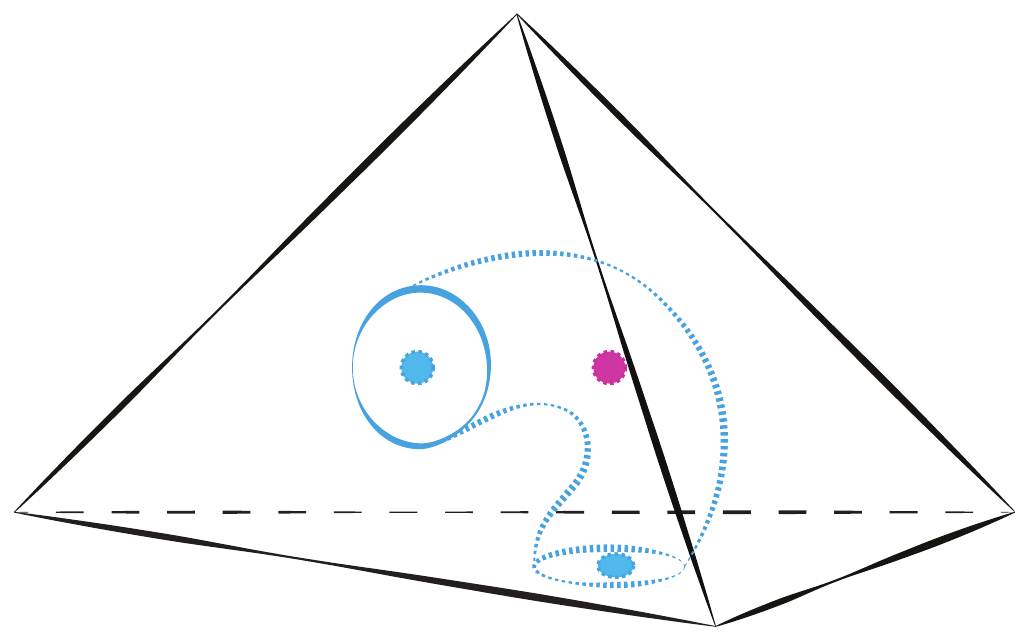}      
   \caption{In $\Sd(\Delta_3)$, $V_\epsilon$ can be a  disc or a cylinder.}
    \label{Complexity2}
      \end{center}
 \end{figure}

\begin{figure}[h]
   \begin{center}
    \includegraphics[scale=0.4]{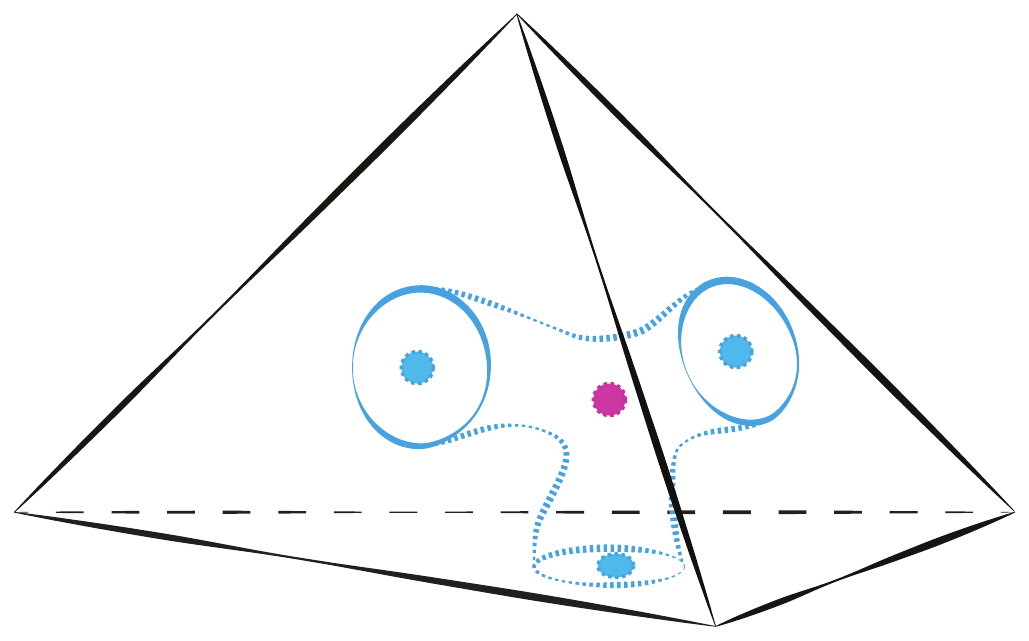}
      \includegraphics[scale=0.4]{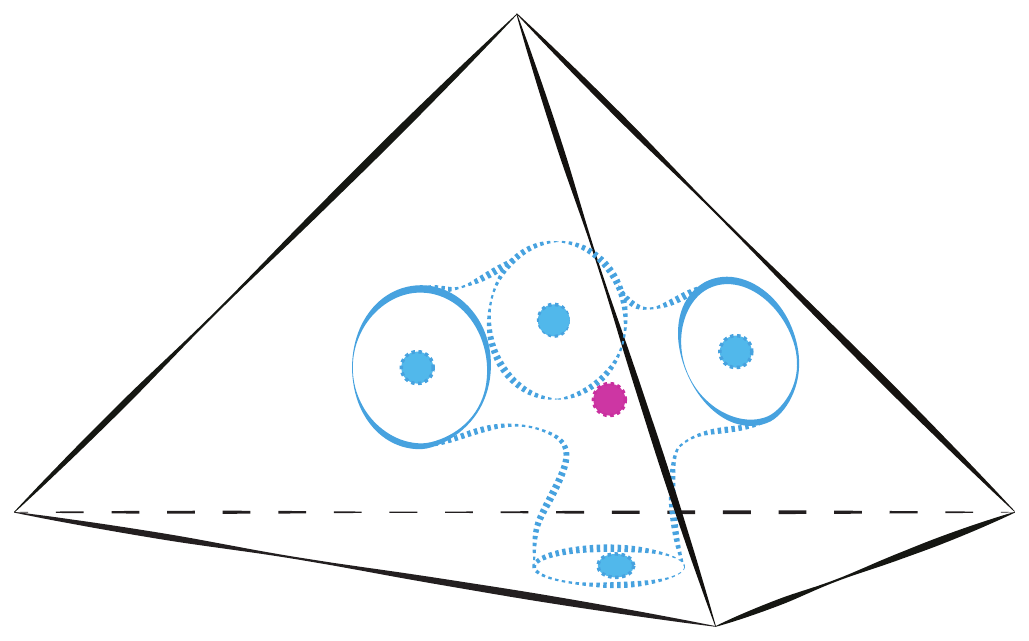}
        \caption{In $\Sd(\Delta_3)$, $V_\epsilon$ a sphere with 3 or 4 holes.}
    \label{Complexity3}
      \end{center}
 \end{figure}

More generally,

\begin{lem}\label{Lem_Complexity} For every $r\in \{0,1,2,3,4\}$,
there exists $\epsilon\in \mC^0(\Sd(\Delta_3))$ such that $V_\epsilon$ is homeomorphic to a sphere with $r$ holes.
\end{lem}
\begin{pf} Let $\epsilon$ take the value 1 on each vertex of $\Delta_3$ and each  barycenter of an edge of $\Delta_3$
and let $\epsilon$ take  the value 0 on the barycenter of $\Delta_3$ itself. Now, depending on whether  $\epsilon$ takes the value 0 or 1 on each barycenter of the codimension-1 faces  of $\Delta_3$, $V_\epsilon$ becomes homeomorphic  to a sphere with up to four holes, see Figures~\ref{Complexity}, \ref{Complexity2}, and \ref{Complexity3}. 
\end{pf}

\begin{thm}\label{Thm_Complexity}
Let  $\Sigma$  be a compact connected orientable surface of Euler characteristic $\chi(\Sigma)\geq 4(3!)^{d-1}-2(4!)^{d-1}$  with $d\geq 1$. Then, there exists $\epsilon\in \mC^0(\Sd^d(\Delta_3))$ such that $V_\epsilon$ is homeomorphic to $\Sigma$.  
\end{thm}
For example, a compact connected orientable surface of genus $0<g\leq 13$ (respectively $13<g\leq 505$) has embeddings of complexity two (respectively three) in $\R^3$, in the sense of Definition~\ref{Defn_Complexity}.

\begin{pf}
We proceed as in Lemma~\ref{Lem_Complexity}. Let $\epsilon$ take the value  1 on each vertex and on the barycenter of each edge of $\Sd^{d-1}(\Delta_3)$ and take the value 0 on the barycenter of each 3-simplex of $\Sd^{d-1}(\Delta_3)$. The number of such 3-simplices is $(4!)^{d-1}.$ If we let $\epsilon$ be 1 on the barycenter of each codimension-2  face of $\Sd^{d-1}(\Delta_3)$, then $V_\epsilon$ becomes homeomorphic to the disjoint union of $(4!)^{d-1}$ copies of the 2-sphere. Changing this value
 to 0 on the barycenter of one interior triangle of $\Sd^{d-1}(\Delta_3)$ results in a connected sum of the two corresponding spheres, which gives rise to a decrease in the Euler characteristic by two. Now the number of 2-dimensional  faces of $\Sd^{d-1}(\Delta_3)$ which lie in the interior of $\Delta_3$ is $\frac{1}{2}(4(4!)^{d-1}-4(3!)^{d-1})=2((4!)^{d-1}-(3!)^{d-1})$, since every such face bounds two 3-simplices and each  such 3-simplex has 4 codimension 1 faces. By letting $\epsilon$ be 1 on the barycenter of each 2-simplex on the boundary of $\Sd(\Delta_3)$ and 0 or 1  on the barycenter of interior ones, we may thus connect sum together the disjoint union of $(4!)^{d-1}$ copies of $S^2$ using up to $2((4!)^{d-1}-(3!)^{d-1})$ cylinders. The first $(4!)^{d}-1$ connected sums can be made to connect together the $(4!)^{d-1}$ copies of  $S^2$ to get a single $S^2$. The result follows.
 \end{pf}

\begin{rem}
 The proof can be carried out in higher dimensions as well, to produce hypersurfaces which are connected  sums of spheres $S^{n-1}$ with handles $S^1\times S^{n-2}$.
\end{rem}

\appendix

\section{More on the universal constants $c_i^+(n,k)$} \label{Ap_ci}
\subsection{Section~\ref{Ssect_$k>1$} revisited}\label{Aps_rev}
For every $1\leq k\leq p$, we set

\begin{equation*} 
\delta^{p,k}(T)=\E_{\nu}(q_{V_\epsilon\cap \stackrel{\circ}{\Delta}_p}(T))=\sum\limits_{i=0}^{p-k}\delta_i^{p,k} T^i,
\end{equation*}
where $\delta^{p,k}_{i}=\int\limits_{\epsilon\in \mC^{k-1}(\Delta_p)}f_i(V_\epsilon\cap \stackrel{\circ}{\Delta}_p)d\mu_\nu(\epsilon).$

These universal polynomials are associated to the standard simplices $(\Delta_p)_{p\geq 1}.$
Recall that for every $\epsilon \in \mC^{k-1}(\Delta_p),$ $V_\epsilon$ is a subcomplex of  $\Sd(\Delta_p)$ so that by 
$V_\epsilon\cap \stackrel{\circ}{\Delta}_p$, we mean the simplices of $V_\epsilon$, the interior of which lie in $\stackrel{\circ}{\Delta}_p$. These are the simplices of the form $[\hat{\sigma}_0,\ldots, \hat{\sigma}_i]$ where $\sigma_i=\Delta_p$.

The first part of Theorem~\ref{Thm_Enu} can be formulated in terms of those polynomials as follows.

\begin{thm} \label{Thm_Enu2} Let $K$ be a finite $n$-dimensional simplicial complex and $k\in\{ 1,\ldots, n\}$. Then, for every $\nu\in[0,1]$, $$\E_{\nu}(q_{V_\epsilon}(T))=\sum_{p=k}^{n}f_p(K)\delta^{p,k}(T).$$
\end{thm}

In order the prove Theorem~\ref{Thm_Enu2}, we need first the following lemma.

\begin{lem}\label{Lem_$E(f_i)$}
Let $K$ be a finite $n$-dimensional simplicial complex and  $L$ be a union of simplices of $K$. Then, for every $\nu\in [0,1]$, every $0\leq i\leq n-k$ and every $k\in \{1,\ldots, n\}$,

$$\E_\nu(f_i(V_\epsilon\cap \stackrel{\circ}{L}))=\sum\limits_{\sigma\in \Sd(K)^{[i]}\cap \stackrel{\circ}{L}} \mu_\nu\{\epsilon\in \mC^{k-1}(K) | \sigma \in V_\epsilon\}. $$
\end{lem}

\begin{pf} We observe that
$$
\begin{array}{lcl}
\E_\nu(f_i(V_\epsilon \cap\stackrel{\circ}{L}))&=&\int_{\mC^{k-1}(K)} f_i (V_\epsilon\cap \stackrel{\circ}{L}) d\mu_\nu(\epsilon)\\
&&\\
&=&\sum\limits_{\epsilon\in \mC^{k-1}(K)}\sum\limits_{\tiny{\begin{array}{c}\sigma\in {\Sd(K)}^{[i]}\cap \stackrel{\circ}{L}\\ \sigma\in V_\epsilon\end{array}}} \mu_\nu(\epsilon)\\
&&\\
&=& \sum\limits_{\sigma\in \Sd(K)^{[i]}\cap \stackrel{\circ}{L}} \sum\limits_{\tiny{\begin{array}{c} \epsilon\in \mC^{k-1}(K)\\\sigma\in V_\epsilon\end{array}}} \mu_\nu(\epsilon)\\
&&\\
&=&\sum\limits_{\sigma\in \Sd(K)^{[i]}\cap \stackrel{\circ}{L}} \mu_\nu\{\epsilon\in \mC^{k-1}(K) | \sigma \in V_\epsilon\},
\end{array}
$$
\end{pf}

In case $L=K$,  Lemma~\ref{Lem_$E(f_i)$} gives
$$\E_\nu(f_i)=\sum\limits_{\sigma\in \Sd(K)^{[i]}} \mu_\nu\{\epsilon\in \mC^{k-1}(K) | \sigma \in V_\epsilon\}. $$

\textbf{Proof of Theorem~\ref{Thm_Enu2}.}
From Lemma~\ref{Lem_$E(f_i)$}, we know that for every $0\leq i \leq  n-k,$

$$\E_{\nu}(f_i)=\sum_{\sigma\in \Sd(K)^{[i]}}\mu_\nu\{\epsilon\in \mC^{k-1}(K)| \sigma\in V_\epsilon\}.$$

Recall that the $i$-simplex $\sigma =[\hat{\sigma}_0,\ldots,\hat{\sigma}_i]\in \Sd(K)$ belongs to $V_\epsilon$ if and only if $\dim \sigma_0\geq k$ and there exists a $k$-face $\tau$ of $\sigma_0$ such that 
$<d\epsilon, \tau>\neq 0.$

In particular, we deduce that $p=\dim \sigma_i\geq k+i$ and

\begin{align}
\E_{\nu}(f_i)&=\sum\limits_{p={k+i}}^n\sum\limits_{\sigma_i\in K^{[p]}}\sum\limits_{\sigma \in \Sd(\sigma_i)^{[i]}\cap \stackrel{\circ}{\sigma_i}}\mu_\nu\{\epsilon\in \mC^{k-1}(K)| \sigma\in V_\epsilon\} \nonumber\\
&= \sum\limits_{p=k+i}^{n}f_p(K)\delta^{p,k}_{i} \label{Eqn_Efi}
\end{align}
since from Lemma~\ref{Lem_$E(f_i)$} it follows that the $i^{th}$ coefficient $\delta_i^{p,k}$ of $\delta^{p,k}(T)$ equals 
$\sum_{\sigma\in \Sd(\Delta_p)^{[i]} \cap \stackrel{\circ}{\Delta}_p}\mu_\nu\{\epsilon \in \mC^{k-1}(\Delta_p)| \sigma\in V_\epsilon\}$, while $\mu_\nu$ is a product measure.

Now, 
$$
\begin{array}{lcl}
\E_{\nu}(q_{V_\epsilon}(T))&=&\sum\limits_{i=0}^{n-k}\E_{\nu}(f_i)T^i\\
&&\\
&=& \sum\limits_{i=0}^{n-k} T^i \sum\limits_{p=k+i}^{n}f_p(K)\delta^{p,k}_i\\
&&\\
&=& \sum\limits_{i=0}^{n-k} T^i \sum\limits_{p=i}^{n-k}f_{p+k}(K)\delta^{p+k,k}_i\\
&&\\
&=& \sum\limits_{p=0}^{n-k} f_{p+k}(K) \sum\limits_{i=0}^{p}\delta^{p+k,k}_iT^i\\
&&\\
&=& \sum\limits_{p=k}^{n} f_{p}(K) \delta^{p,k}(T).\\
\end{array}
$$

\hfill$\square$

Let us set, for every $\nu\in [0,1]$ and $k\leq p\leq n$, $\tau_p=\delta^{p,k}(-1).$

\begin{cor} \label{Cor_Enuch}Let $K$ be a finite $n$-dimensional simplicial complex and $k\in\{1,\ldots, n\}$, $\nu\in[0,1]$. Then, $$\E_{\nu}(\chi)=\sum_{p=k}^nf_p(K)\tau_p.$$

Moreover, for every $0\leq i\leq n-k$, the following average Morse inequalities hold
$$\E_{\nu}(b_i)\leq \E_{\nu}(f_i)=\sum_{p=k+i}^nf_p(K)\delta^{p,k}_i$$ and
$$\sum_{j=0}^i(-1)^{i-j}\E_{\nu}(b_j)\leq \sum_{j=0}^i(-1)^{i-j}\E_{\nu}(f_j),$$ where the latter is an equality if $i=n-k$.
\end{cor}

From Corollary~\ref{Cor_Enuch} we deduce another formulation of the last part of Corollary~\ref{Cor_EbiEuler$k>1$}, since  we deduce that $\sum_{i=0}^{n-k}(-1)^i c_i^+(n,k)=\sum_{p=k}^{n}\tau_pq_{p,n}.$

\begin{pf} The first part follows from Theorem~\ref{Thm_Enu2} after evaluation at $T=-1$, since by definition, for every $\epsilon\in \mC^{k-1}(K),$ $\chi(V_\epsilon)=q_{V_\epsilon}(-1)$. Then, for every $\epsilon\in \mC^{k-1}(K)$, the Morse inequalities applied to the simplicial  chain complex of $V_\epsilon$ with $\Z/2\Z$-coefficients read $b_i(V_\epsilon)\leq f_i(V_\epsilon)$ and $\sum_{j=0}^i (-1)^{i-j}b_j(V_\epsilon)\leq \sum_{j=0}^i(-1)^{i-j}f_j(V_\epsilon)$, the latter being an equality when $i=n-k$. The last part of Corollary~\ref{Cor_Enuch} thus follows from (\ref{Eqn_Efi}) after integration over $ \mC^{k-1}(K).$
\end{pf}

\textbf{Examples :}
\begin{enumerate}
\item When $k=1$ the first part of Corollary~\ref{Cor_Enuch} gives back Corollary~\ref{Cor_E(chi)1}, as follows from Corollary~\ref{Cor_tau}. 
\item When $k=n$ and $\nu=\frac{1}{2}$, $V_\epsilon$ is a finite set of points for every $\epsilon \in \mC^{n-1}(K)$ and Corollary~\ref{Cor_Enuch} combined with  Corollary~\ref{Cor_tau} gives $\E(\chi)=\frac{1}{2}f_n(K)$.
\item When $k=n-1$ and $\nu=\frac{1}{2}$, $V_\epsilon$ is a graph for every $\epsilon \in \mC^{n-2}(K)$ and the second part of Corollary~\ref{Cor_Enuch} combined with Corollary~\ref{Cor_tau} gives

$\begin{array}{lcl}
\E(b_0)&\leq & \frac{f_{n-1}(K)}{2}+(1-\frac{1}{2^n})f_n(K),\\
&&\\
\E(b_1)&\leq & \frac{n+1}{2}f_n(K),\\
&&\\
\E(\chi)&= & \frac{f_{n-1}(K)}{2}+(\frac{1-n}{2}-\frac{1}{2^n})f_n(K).
\end{array} $
\end{enumerate}

\subsection{Computations of the universal polynomials  $\delta^{p,k}(T)$}\label{deltapk}

The coefficients of the universal polynomials $\delta^{p,k}(T)$ introduced in Section~\ref{Aps_rev}  are given by the following theorem.

\begin{thm} \label{Thm_deltatau}
For  every $k\in\{1,\ldots, n\}$, every $i\in \{0,\ldots, p-k\}$ and every $\nu\in[0,1]$, 
$$ \delta^{p,k}_i=\sum\limits_{j=0}^{i}\binom{i}{j}(-1)^{i-j} \sum\limits_{l=i}^{p-k}\binom{p+1}{l}\big(1-\mu_\nu(Z^{k-1}(\Delta_{p-l}))\big)j^{l}$$ 

and

$$\tau_p=\sum\limits_{l=k}^{p}\binom{p+1}{l+1}(-1)^{p-l} \big(1-\mu_\nu(Z^{k-1}(\Delta_l))\big).$$ 
\end{thm}

In particular, we deduce from  Corollary~\ref{Cor_MuZ} and Theorem~\ref{Thm_deltatau} that for $\nu=\frac{1}{2}$, $$\delta^{p,k}_i=\sum_{j=0}^{i}\binom{i}{j}(-1)^{i-j}\sum\limits_{l=i}^{p-k}\binom{p+1}{l}\big(1-\frac{1}{2^{\binom{p-l}{k}}}\big)j^{l}$$
and $$\tau_p=\sum_{l=k}^{p}\binom{p+1}{l+1}(-1)^{p-l}\big(1-\frac{1}{2^{\binom{l}{k}}}\big).$$

 In order to prove Theorem~\ref{Thm_deltatau}, we need first the following Lemma~\ref{Lem_Indice}.

For every $l<p$ and every $0< i\leq p-l$, let us denote by $\ind(l,p,i)$ the number of $i$-simplices of $\Sd(\Delta_p)$ which are of the form $[\hat{\sigma}_0,\ldots,\hat{\sigma}_i]$ with $\dim\sigma_0=l$ and $\dim \sigma_i=p$.

\begin{lem} \label{Lem_Indice}
For every $l<p$ and every $0< i\leq p-l$, $\ind(l,p,i)=\binom{p+1}{l+1}\lambda_{p-l,i}$.
\end{lem}

Recall that  $\lambda_{p-l,i}$  is the number of  interior $(i-1)$-faces 
of $\Sd(\Delta_{p-l-1})$.

\begin{pf} We observe that 
\begin{enumerate}
\item there are $\binom{p+1}{l+1}$ choices for   a $l$-simplex $\sigma_0$ of $\Delta_p$.  
\item there is a bijection between the $(i+1)$-flags ($\sigma_0<\ldots<\sigma_i=\Delta_p$) and the $i$-flags  ($\Lk({\sigma}_0, \sigma_1)<\dots<\Lk({\sigma}_0,\sigma_i)=\Delta_{p-l-1})$.
\end{enumerate}

Now, the $i$-flags ($\Lk({\sigma}_0, \sigma_1)<\dots<\Lk({\sigma}_0,\sigma_i)$)
 exactly define the $(i-1)$-simplices interior to $\Sd(\Delta_{p-l-1})$. By definition there are $\lambda_{p-l,i}$ many such simplices, hence the result. \end{pf}

\textbf{Proof of Theorem~\ref{Thm_deltatau}.} By definition and Lemma~\ref{Lem_$E(f_i)$},
$$\delta^{p,k}_i=\sum_{\sigma\in\Sd(\Delta_p)^{[i]}\cap \sDp} \mu_\nu\{\epsilon\in \mC^{k-1}(\Delta_p)| \sigma\in V_\epsilon\}.$$
If $i=0$, we deduce that $\delta^{p,k}_0=1-\mu_\nu(Z^{k-1}(\Delta_p))$.
If $i>0$, a simplex $\sigma\in \Sd(\Delta_p)^{[i]}\cap \stackrel{\circ}{\Delta_p}$  is of the form $[\hat{\sigma}_0,\ldots,\hat{\sigma}_i]$ where $\sigma_i=\Delta_p$ and $\sigma_0<\sigma_i$ is a face of dimension $l\in\{0,\ldots,p-i\}$. Moreover, such a simplex belongs to $V_\epsilon$ if and only if $l\geq k$ and  the restriction of $d\epsilon$ to $\sigma_0$ does not vanish. We thus deduce

\begin{align}
\delta^{p,k}_i&=\sum\limits_{l=k}^{p-i}\ind(l,p,i)\big(1-\mu_\nu(Z^{k-1}(\Delta_l))\big)\nonumber\\
&\nonumber\\
&=\sum\limits_{l=k}^{p-i} \binom{p+1}{l+1}\lambda_{p-l,i} \big(1-\mu_\nu(Z^{k-1}(\Delta_l))\big) \label{Eqn_delta}\\
&\nonumber\\
&=\sum\limits_{l=i}^{p-k} \binom{p+1}{l}\lambda_{l,i} \big(1-\mu_\nu(Z^{k-1}(\Delta_{p-l}))\big),\nonumber
\end{align}
where the second line follows from Lemma~\ref{Lem_Indice}.
From  \cite{DPS}, (see also \cite{SW1}) we now deduce,

$$\begin{array}{lcl}
\delta_i^{p,k}&=&\sum\limits_{l=i}^{p-k}\binom{p+1}{l}\sum\limits_{j=0}^{i}\binom{i}{j}(-1)^{i-j}j^l\big(1-\mu_\nu(Z^{k-1}(\Delta_{p-l}))\big)\\
&&\\
&=&\sum\limits_{j=0}^{i}\binom{i}{j}(-1)^{i-j} \sum\limits_{l=i}^{p-k}\binom{p+1}{l}\big(1-\mu_\nu(Z^{k-1}(\Delta_{p-l}))\big)j^{l},\\
\end{array}
$$
where the formula remains valid for $i=0$ with the convention $0^0=1.$
Then,

$\begin{array}{lcl}
\tau_p&=&\sum\limits_{i=0}^{p-k}(-1)^i\delta^{p,k}_i=\sum\limits_{i=0}^{p-k}(-1)^i\sum\limits_{l=i}^{p-k} \binom{p+1}{l}\lambda_{l,i}\big(1-\mu_\nu(Z^{k-1}(\Delta_{p-l}))\big)\\\\
&&\\
&=&\sum\limits_{l=0}^{p-k} \binom{p+1}{l}\big(1-\mu_\nu(Z^{k-1}(\Delta_{p-l}))\big)\sum\limits_{i=0}^{l}(-1)^i\lambda_{l,i}\\
&&\\
&=&1-\mu_\nu(Z^{k-1}(\Delta_p))+\sum\limits_{l=1}^{p-k} \binom{p+1}{l} \big(1-\mu_\nu(Z^{k-1}(\Delta_{p-l}))\big)\sum\limits_{i=0}^{l}(-1)^i\lambda_{l,i}\\
&&\\

\end{array}$

Since by definition $\lambda_{l,0}=0$ for every $l>0$. Moreover,
$\sum_{i=1}^{l}(-1)^i\lambda_{l,i}=-(\chi(\Delta_{l-1})-\chi({\partial\Delta_{l-1}}))=(-1)^l.$ 
We thus deduce

$\begin{array}{lcl}
&&\\
\tau_p&=&\sum\limits_{l=0}^{p-k} \binom{p+1}{l} \big(1-\mu_\nu(Z^{k-1}(\Delta_{p-l}))\big)(-1)^l\\
&&\\
&=&\sum\limits_{l=k}^{p}\binom{p+1}{l+1}(-1)^{p-l} \big(1-\mu_\nu(Z^{k-1}(\Delta_l))\big).\\
\end{array}$

\hfill $\square$

\begin{cor}\label{Cor_tau}

\begin{enumerate}

\item  Let $\nu=\frac{1}{2}$. Then, for every $0\leq i\leq n-k$, $\delta^{k+i,k}_i=\frac{i!}{2}\binom{k+i+1}{i}$.

\item Let $\nu=\frac{1}{2}$. Then, for every $1\leq k \leq p\leq n$, $\delta^{p,k}_0=1-\frac{1}{2^{\binom{p}{k}}}$, and $\tau_k=\frac{1}{2}$.

\item If $k=1$, $\tau_p=(-1)^{p+1}(1-\nu^{p+1}-(1-\nu)^{p+1})$ for every $1\leq p \leq n$ and every $\nu\in[0,1]$.  

 \end{enumerate}
\end{cor}

\begin{pf} 

\begin{enumerate}

\item From Corollary~\ref{Cor_MuZ} and Theorem~\ref{Thm_deltatau},  $\delta^{k+i,k}_i=\sum\limits_{l=i}^{i}\binom{k+i+1}{l} \lambda_{l,i} (1-\frac{1}{2^{\binom{k+i-l}{k}}})= \binom{k+i+1}{i} \lambda_{i,i} (1-\frac{1}{2})$. The result then follows from 
\cite{SW1}.
\item  Since $\lambda_{l, 0}=0$ for every $l>0$ and $\lambda_{0,0}=1$, Corollary~\ref{Cor_MuZ} and Theorem~\ref{Thm_deltatau} imply that $\delta^{p,k}_0= 1-\frac{1}{2^{\binom{p}{k}}},$  and   $\tau_k=(1-\frac{1}{2})=\frac{1}{2}$.
\item  When $k=1$, we know from  Theorem~\ref{Thm_deltatau} that for every $\nu\in[0,1]$, 
$$\begin{array}{lcl}
\tau_p&=&\sum\limits_{l=0}^{p-1}(-1)^l\binom{p+1}{l}(1-\nu^{p-l+1}-(1-\nu)^{p-l+1}).
\end{array}$$

Thus, $\tau_p$ is written as the sum of three terms $A_1=\sum_{l=0}^{p-1}(-1)^l\binom{p+1}{l}$, $A_2=-\nu^{p+1}\sum_{l=0}^{p-1}(\frac{-1}{\nu})^l\binom{p+1}{l}$ and $A_3=-(1-\nu)^{p+1}\sum_{l=0}^{p-1}(\frac{-1}{1-\nu})^l\binom{p+1}{l}$.

Using binomial expansion we get

$\begin{array}{l}
A_1=(1-1)^{p+1}-(-1)^pp-(-1)^{p+1},\\
A_2=-\nu^{p+1}((1-\frac{1}{\nu})^{p+1}-(\frac{-1}{\nu})^pp-(\frac{-1}{\nu})^{p+1}),\\
A_3=-(1-\nu)^{p+1}((1-\frac{1}{1-\nu})^{p+1}-(\frac{-1}{1-\nu})^pp-(\frac{-1}{1-\nu})^{p+1}).
\end{array}$

Hence the result.
\end{enumerate}
\end{pf}

\bibliography{Random2}
\bibliographystyle{abbrv}

Univ Lyon, Universit\'e Claude Bernard Lyon 1, CNRS UMR 5208, Institut Camille Jordan, 43 blvd. du 11 novembre 1918, F-69622 Villeurbanne cedex, France

{salepci@math.univ-lyon1.fr, welschinger@math.univ-lyon1.fr.}
\end{document}